\documentclass[11pt]{amsart}
\usepackage{amsmath,amsxtra,amssymb,amsthm,amsfonts}
\usepackage{color}

\numberwithin{equation}{section}
\newtheorem{lemma}{Lemma}[section]
\newtheorem{prop}[lemma]{Proposition}
\newtheorem{theorem}[lemma]{Theorem}

\newtheorem{cor}[lemma]{Corollary}

\newtheorem{rem}[lemma]{Remark}

\newcommand{\re}{\begin{rem}\rm}
  \newcommand{\mar}{\end{rem}}

\newtheorem{defi}[lemma]{Definition}

\newcommand{\kla}{\left ( }
\newcommand{\mer}{\right ) }

\renewcommand{\for}{\begin{eqnarray*}}

\newcommand{\mel}{\end{eqnarray*}}

\newcommand{\kl}{\pl \le \pl}
\newcommand{\gl}{\pl \ge \pl}

\newcommand{\lel}{\pl = \pl}

\newcommand{\ran}{\rangle}
\newcommand{\lan}{\langle}

\newcommand{\nz}{{\mathbb N}}
\newcommand{\nen}{n \in \nz}

\newcommand{\rz}{{\mathbb R}}
\newcommand{\zz}{{\mathbb Z}}

\newcommand{\cz}{{\mathbb C}}

\newcommand{\ten}{\otimes}

\newcommand{\pl}{\hspace{.1cm}}

\newcommand{\qd}{\end{proof}\vspace{0.5ex}}

\newcommand{\Om}{\Omega}
\newcommand{\om}{\omega}

\newcommand{\al}{\alpha}

\newcommand{\si}{\sigma}

\newcommand{\la}{\lambda}
\newcommand{\eps}{\varepsilon}

\newcommand{\F}{{\mathcal F}}

\newcommand{\A}{{\mathcal A}}

\newcommand{\R}{{\mathcal R}}

\newcommand{\U}{{\mathcal U}}

\newcommand{\pf}{\begin{proof}}

\newcommand{\xspace}{\hbox{\kern-2.5pt}}
\newcommand{\xyspace}{\hbox{\kern-1.1pt}}

\newcommand{\TT}{{\mathbb T}}

\newcommand{\ZZ}{\mathbb{Z}}
\newcommand{\RR}{\mathbb{R}}

\allowdisplaybreaks

\definecolor{LightGray}{rgb}{0.94,0.94,0.94}
\definecolor{VeryLightBlue}{rgb}{0.9,0.9,1}
\definecolor{LightBlue}{rgb}{0.8,0.8,1}
\definecolor{DarkBlue}{rgb}{0,0,0.6}
\definecolor{LightGreen}{rgb}{0.88,1,0.88}
\definecolor{MidGreen}{rgb}{0.6,1,0.6}
\definecolor{DarkGreen}{rgb}{0,0.6,0}
\definecolor{DarkGrreen}{rgb}{0,0.8,0}
\definecolor{VeryLightYellow}{rgb}{1,1,0.9}
\definecolor{LightYellow}{rgb}{1,1,0.6}
\definecolor{MidYellow}{rgb}{1,1,0.5}
\definecolor{DarkYellow}{rgb}{0.8,1,0.3}
\definecolor{VeryLightRed}{rgb}{1,0.9,0.9}
\definecolor{LightRed}{rgb}{1,0.8,0.8}
\definecolor{DarkRed}{rgb}{0.8,0.2,0}
\definecolor{DarkRedb}{rgb}{0.6,0.2,0}
\definecolor{DarkLila}{rgb}{0.8,0,1}
\definecolor{Beige}{rgb}{0.96,0.96,0.86}
\definecolor{Gold}{rgb}{1.,0.84,0.}
\definecolor{Goldb}{rgb}{0.7,0.3,0.5}
\definecolor{MyYellow}{rgb}{1.,0.84,0.8}

\begin{document}

\title[]{Some Classification Results for Generalized $q$-Gaussian Algebras}

\author[Marius Junge]{Marius Junge$^*$}
\address{Department of Mathematics\\
University of Illinois, Urbana, IL 61801, USA} \email[Marius Junge]{mjunge@illinois.edu}

\author[Stephen Longfield]{Stephen Longfield}
\address{Department of Mathematics\\
University of Illinois, Urbana, IL 61801, USA} \email[Stephen Longfield] {longfie2@illinois.edu}

\author[Bogdan Udrea]{Bogdan Udrea}
\address{Department of Mathematics\\
University of Illinois, Urbana, IL 61801, USA\\ \newline 
\indent and  IMAR, Bucharest, Romania} \email[Bogdan Udrea]{budrea@illinois.edu}

\thanks{$^*$ Marius Junge is partially supported by  NSF-DMS 1201886}

\begin{abstract} To any trace preserving action $\si: G \curvearrowright A$ of a countable discrete group on a finite von Neumann algebra $A$ and any orthogonal representation $\pi:G \to \mathcal O(\ell^2_{\rz}(G))$, we associate the generalized q-gaussian von Neumann algebra $A \rtimes_{\si} \Gamma_q^{\pi}(G,K)$, where $K$ is an infinite dimensional separable Hilbert space. Specializing to the cases of $\pi$ being trivial or given by conjugation, we then prove that if $G \curvearrowright A = L^{\infty}(X)$, $G' \curvearrowright B = L^{\infty}(Y)$ are p.m.p. free ergodic rigid actions, the commutator subgroups $[G,G]$, $[G',G']$ are ICC, and $G, G'$ belong to a fairly large class of groups (including all non-amenable groups having the Haagerup property), then $A \rtimes \Gamma_q(G,K) = B \rtimes \Gamma_q(G',K')$ implies that $\mathcal R(G \curvearrowright A)$ is stably isomorphic to $\mathcal R(G' \curvearrowright B)$, where $\mathcal R(G \curvearrowright A), \mathcal R(G' \curvearrowright B)$ are the countable, p.m.p. equivalence relations implemented by the actions of $G$ and $G'$ on $A$ and $B$, respectively. Using results of D. Gaboriau and S. Popa we construct continuously many pair-wise non-isomorphic von Neumann algebras of the form $L^{\infty}(X) \rtimes \Gamma_q(\mathbb{F}_n,K)$, for suitable free ergodic rigid p.m.p. actions $\mathbb{F}_n \curvearrowright X$.
\end{abstract}

\maketitle\tableofcontents

\section{Introduction}Ever since Murray and von Neumann laid the foundations of von Neumann algebras theory, classification of the objects involved (then called rings of operators) was a crucial issue. Specifically, the first non-trivial examples of factors were of the form $L^{\infty}(X)\rtimes \Gamma$ (the so-called group measure space construction) for certain actions of discrete countable groups on measure spaces or of the form $L(\Gamma)$ (group von Neumann algebras) for countable, discrete, ICC groups.  In this paper we want to study and eventually classify similar von Neumann algebras constructed by data from a group action and $q$-Gaussian algebra.

Our results are motivated by the success in classifying group von Neumann algebras through Popa's deformation-rigidity theory, which we will now review. In fact, the first natural question is whether the crossed products are completely classified by their original data, i.e. if isomorphism of two cross-product factors implies isomorphism of the original actions, or if isomorphism of group factors forces the groups to be isomorphic. When this ambitious goal is not attainable, one might still ask if isomorphism of the associated von Neumann algebras forces the actions or the groups to share some properties, even without being isomorphic. Two p.m.p. actions $\Gamma \curvearrowright X$, $\Lambda \curvearrowright Y$ are called isomorphic (or \emph{conjugate}) if there exist a measure space isomorphism $\Delta : X \to Y$ and a group isomorphism $\delta: \Gamma \to \Lambda$ such that $\Delta(gx)=\delta(g)\Delta(x)$, for every $g \in \Gamma$ and almost every $x \in X$. A weaker equivalence between two actions is that of \emph{orbit equivalence} (or OE): $\Gamma \curvearrowright X$ and $\Lambda \curvearrowright Y$ are called OE if there exists a measure space isomorphism $\Delta:X \to Y$ such that $\Delta(\Gamma x)=\Lambda \Delta(x)$ for almost every $x \in X$. This can be relaxed to \emph{stable orbit equivalence}, which means the existence of Borel subsets $Z \subset X, T \subset Y$ that intersect alomost every orbit and of a nonsingular isomorphism $\Delta:Z \to T$ such that $\Delta(\Gamma x \cap Z)=\Lambda \Delta(x) \cap T$, for almost every $x \in Z$. A still weaker notion is that of von Neumann equivalence (VNE): $\Gamma \curvearrowright X$ and $\Lambda \curvearrowright Y$ are called von Neumann equivalent if $L^{\infty}(X)\rtimes \Gamma \cong L^{\infty}(Y)\rtimes \Lambda$. It should be noted that by a result of Singer (\cite{Sing55}), OE amounts to the existence of an isomorphism between the two cross-product von Neumann algebras carrying $L^{\infty}(X)$ onto $L^{\infty}(Y)$. Using this precise terminology, the complete classification goal amounts to being able to prove that if two actions are VNE then they are conjugate. A weaker result would be obtained by proving that if they are OE, then they have to be conjugate.
\par It became gradually clear over a long period of time that none of the implications above holds in full generality. This culminated in Connes' ground-breaking result that all type $II_1$ injective factors are isomorphic (\cite{Connes76}), which leads to the conclusion that any p.m.p., free, ergodic action of any countable discrete amenable group gives rise to one and the same von Neumann algebra - the type $II_1$ hyperfinite factor, and also that for any ICC, countable discrete amenable group $\Gamma$, $L(\Gamma)$ is isomorphic to the same hyperfinite factor. On the other hand, Connes, Feldman and Weiss proved that all (free, ergodic, p.m.p.) actions of (discrete, countable) amenable groups are OE (\cite{CFW}). So within the realm of cross-product factors arising from actions of amenable groups, one cannot distinguish the objects at all in terms of their original data. In colloquial terms, a factor $L^{\infty}(X)\rtimes \Gamma$ with $\Gamma$ amenable ``remembers" nothing about the action or the group, except that the group is amenable. The classification goal prompts one to look for a ``rigidity" phenomenon (broadly speaking), i.e. when the von Neumann algebra remembers some amount of information (ideally everything) about its building data. To make this more precise, an action is called
\begin{itemize}
\item OE-superrigid if any other action which is OE to it must be conjugate to it;
\item W*-superrigid if any other action which is VNE to it must be conjugate to it.
\end{itemize}
Thus, the associated von Neumann algebra completely remembers the action in the case of W*-superrigid actions. In the case of OE-superrigid actions, the action can be reconstructed from its OE class. Along these lines, let us mention here \cite{Zimmer, Fur99, PoCSR, PoSG, IoaCSR, IoaWS, Ki10, Ki11}.
\par Some fifteen years ago, Popa's deformation rigidity theory began to produce the first significant results in this direction (\cite{PoBe}, \cite{PoI}, \cite{PoCSR}, \cite{PoSG}). Since then  ground-breaking results have been obtained by Popa and his collaborators, see e.g. \cite{IPP}, \cite{IoaWS}, \cite{IPV}, \cite{ioana2011}, \cite{OPCartanI}, \cite{OPCartanII}, \cite{CS}, \cite{PetL2}, \cite{PeWS}, \cite{PoVaI}. To cite only a few, Popa proved strong rigidity result for cross-product factors which come from Bernoulli actions of w-rigid groups (see \cite{PoI,PoII}). Then, in \cite{PoCSR, PoSG} he proved cocycle superrigidity results for malleable actions (notably Bernoulli) of either rigid groups or having the spectral gap property (e.g. for direct products $H\times G$ with $H$ infinite and $G$ non-amenable). This in particular implies that any (free ergodic) action which is OE to a Bernoulli action of such a group has to be conjugate to it. This was further upgraded by Ioana, who proved in \cite{IoaWS} that the Bernoulli actions of property (T) groups are (virtually) W*-superrigid. In the same vein, Popa and Vaes found the first examples of groups for which every action gives rise to a factor having unique group measure space Cartan subalgebra, which, when coupled with Kida's OE superrigidity results in \cite{Ki11} lead to the first examples of groups whose every action is W*-superrigid (\cite{PoVaWS}), results further extended in \cite{CP}. On the other hand, the ground-breaking results of Ozawa and Popa in \cite{OPCartanI,OPCartanII} provided the first examples of type $II_1$ factors having unique Cartan subalgebra, results further enhanced by Chifan and Sinclair in \cite{CS}, then Chifan, Sinclair and the last author in \cite{CSU} and ultimately by Popa and Vaes, who proved in \cite{PoVaI,PoVaII} that every action of any weakly amenable group with positive first Betti number, as well as of any non-amenable hyperbolic group, gives rise to a von Neumann algebra having unique Cartan subalgebra.
\par The q-Gaussian algebras ($-1<q<1$) were defined by Bo\.{z}ejko and Speicher (\cite{BoSpe}, \cite{BoKuSpe}) and studied further by Krolak (\cite{Kro}), Ricard (\cite{Ri}) who proved the factoriality of these algebras and Nou who proved they are non-amenable (\cite{Nou}).  Shlyakhtenko proved solidity of $\Gamma_q(\mathcal{H})$ for $q\le \sqrt{2}-1$ following Ozawa's approach in \cite{Shl1} and absence of Cartan subalgebras for small $q$ introducing the power series approach in \cite{Shl}, see also \cite{GS}. Avsec (\cite{Av}) proved that they have the complete metric approximation property and as byproduct that for $|q|<1$ these algebras are \emph{strongly solid} if dim$\mathcal H < \infty$, (for a definition, see \cite{OPCartanI}) using deformation-rigidity techniques. These algebras  can be thought of both as von Neumann algebra implementations of the canonical q-commutation relations or as interpolations between the classical commutative gaussian random variables (when $q=1$) and the hyperfinite type $II_1$ factor ($q=-1$), going through the free group factors ($q=0$). The $q$-gaussian algebras are probabilistic in nature. Indeed, with the help of the $q$-gaussian relations one can find a distinct family of brownian motions which all satisfy L\'evy's axioms for the classical brownian motions (except for commutativity). In our rigidity related context, the $q$-gaussian relations become gradually more difficult the more commutative they are, or alternatively less free. For $q=0$, the von Neumann algebras $\Gamma_q(\mathcal H)=L(\mathbb{F}_{dim \mathcal H})$ appear in the work of Voiculescu, Dykema and Nica \cite{VDN}.

\par In this paper we prove a ``weak rigidity'' result for certain classes of generalized q-Gaussian algebras with action. Our algebras are a mix of classical q-Gaussian algebras and the cross-product construction  and that is why we use the  suggestive notation $A \rtimes \Gamma_q(G,K)$ for them. To be more precise, for each trace preserving action $\si:G \curvearrowright A$ on a finite von Neumann algebra $A$, every orthogonal representation $\pi:G \to \mathcal{O}(\ell^2_{\rz}(G))$ and every infinite dimensional separable Hilbert space $K$ we construct a \emph{generalized q-Gaussian von Neumann algebra} $A \rtimes_{\si} \Gamma_q^{\pi}(G,K)$ as a suitable subalgebra of a crossed product (see Section 4 for a precise definition). This construction also makes sense in the case of unitary representations on complex Hilbert spaces, but we are dealing mostly with the real case, except for some of our examples in Section 7. The main result we prove is (Theorem 7.2):
\begin{theorem}\label{00}Let $M = A \rtimes \Gamma_q(G,K) = B \rtimes \Gamma_q(G',K')$ with the representation $\pi:G \to \mathcal O (\ell^2_{\rz}(G))$ either trivial or given by conjugation and assume that $A$ and $B$ are abelian, the inclusions $A \subset M$ and $B \subset M$ are rigid, $[G,G]$, $[G',G']$ are ICC groups, and the actions $G \curvearrowright A$, $G' \curvearrowright B$ are free and ergodic. If moreover one of the following conditions holds:
\begin{enumerate}
\item $q=0$;
\item $G,G'$ are groups with the Haagerup property;
\item $\pi$ is trivial, $[G,G]$ and $[G',G']$ are weakly amenable groups which admit unbounded 1-cocycles into mixing non-amenable representations;
\item $\pi$ is trivial, $[G,G]$, $[G',G']$ are weakly amenable groups which admit proper 1-cocycles into non-amenable representations;
\item $\pi$ is trivial, $[G,G]$, $[G',G']$ are weakly amenable, non-amenable bi-exact groups,
\end{enumerate}
then $\mathcal R(G \curvearrowright A)$ and $\mathcal R(G' \curvearrowright B)$ are stably isomorphic.
\end{theorem}

Note that any weakly amenable group (or having the Haagerup property) $G=\Gamma_1 \ast \Gamma_2$, where $|\Gamma_1|\geq 2,|\Gamma_2|\geq 3$ or more generally any weakly amenable non-trivial free product $G=\ast_i G_i$, satisfies the assumptions in item 3, and any $G$ such that $[G,G]$ is a non-amenable hyperbolic group satisfies the assumptions in item 5. Indeed, any non-trivial free product admits an unbounded 1-cocycle into its left regular representation. Now $[G,G]$ is an infinite group (otherwise, since $G/[G,G]$ is abelian, $G$ would follow amenable). The restriction of the cocycle to $[G,G]$ has to remain unbounded, because otherwise by Thm. 2.5 in \cite{CP} the cocycle would be bounded on the whole of $G$, a contradiction. In particular the free groups are good examples in both cases. By the results in \cite{PoGa}, there exist uncountably many stably non-OE free ergodic rigid pmp actions $\mathbb{F}_n \curvearrowright (X,\mu)$. This leads to the following consequence:
\begin{cor}
There exist continuously many pairwise non-isomorphic von Neumann algebras of the form $L^{\infty}(X) \rtimes \Gamma_q(\mathbb{F}_n,K)$.
\end{cor}
Using Thm. 1.3 in \cite{Ga} one can replace the free groups by any weakly amenable (or having the Haagerup property) non-trivial free product group $G=\ast_i G_i$, thereby obtaining
\begin{cor}For any non-trivial free product $G=\ast_i G_i$ which is weakly amenable or has the Haagerup property, there exist continuously many pairwise non-isomorphic type $II_1$ factors of the form $L^{\infty}(X) \rtimes \Gamma_q(G,K)$.
\end{cor}
By exploiting the ``Bass-Serre rigidity'' results in \cite{IPP, CH} we also obtain
\begin{cor}Let $G_1,...,G_m$, $H_1,...,H_n$ be ICC groups, each of which either contains a non-virtually abelian subgroup with relative property (T) or is a direct product of a non-amenable and an infinite group. Denote by $G=G_1 \ast \ldots \ast G_m$, $H=H_1 \ast \ldots \ast H_n$. Assume that $G$ and $H$ are weakly amenable or have the Haagerup property. Let $G \curvearrowright X$, $H \curvearrowright Y$ be two p.m.p. free ergodic rigid actions such that the restriction to each factor is still ergodic. If $L^{\infty}(X)\rtimes \Gamma_q(G,K)$ is isomorphic to $L^{\infty}(Y)\rtimes  \Gamma_q(H,K)$, then $m=n$ and after a permutation of indices we have $\mathcal R(G_i \curvearrowright X)=\mathcal R(H_i \curvearrowright Y)$, for all $i$.
\end{cor}
Using the results of Monod and Shalom in \cite{MoSh}, we also deduce:
\begin{cor}Let $G=\mathbb{F}_{n_1}\times \ldots \times \mathbb{F}_{n_k} \curvearrowright X$, $G'=\mathbb{F}_{m_1}\times \ldots \times \mathbb{F}_{m_l} \curvearrowright X$ be pmp free ergodic rigid actions. If $k\neq l$, then $L^{\infty}(X)\rtimes \Gamma_q(G,K)$ and $L^{\infty}(X)\rtimes \Gamma_q(G',K)$ are non-isomorphic.
\end{cor}
\par Thus if we consider the class $\mathcal{C}_q$ of q-Gaussian von Neumann algebras $A \rtimes \Gamma_q(G,K)$ such that all the conditions in Theorem 1.1 are satisfied, then for two isomorphic objects $M \cong M'$ in $\mathcal{C}_q$, it follows that that the actions $G \curvearrowright A$, $G' \curvearrowright B$ are \emph{stably orbit equivalent} (see e.g. \cite{Fur99}). In particular, if the initial actions are not (stably) OE, then the corresponding generalized q-Gaussians cannot be isomorphic. Our result can be seen as a partial classification result much in the spirit of \cite{PoBe} and \cite{IPP} and the more recent \cite{MeVa}, allowing one to recapture some of the information contained in the original data these von Neumann algebras are built of. Indeed, in \cite{PoBe}, Popa considered the class $\mathcal{HT}_s$ of all type $II_1$ factors $M$ having a Cartan subalgebra $A$ such that the inclusion $A\subset M$ is rigid and $M$ has the Haagerup property relative to $A$. He was able to prove that for two factors $M_{1,2} \in \mathcal{HT}_s$, if $M_1=M_2$, then the corresponding Cartan subalgebras $A_1$ and $A_2$ have to be unitarily conjugate in $M$, and in particular the equivalence relations associated to the inclusions $A_1 \subset M_1$, $A_2 \subset M_2$ are isomorphic. Though we cannot prove that $A$ and $B$ are unitarily conjugate, as they are not MASAs, we are still able to conclude that $\mathcal R(A \subset M) \cong \mathcal R(B \subset M)$, by making crucial use of some recent results of Meesschaert and Vaes (\cite{MeVa}). Here $\mathcal R(A \subset M)$ is the generalized equivalence relation associated to an inclusion $A \subset M$, where $M$ is a type $II_1$ factor and $A$ is an abelian subalgebra which is not maximal abelian (see section 3 and \cite{MeVa}). This generalized equivalence relation does not coincide, in general, with the classical one, when $A$ is not a MASA. However, it turns out that in our case $\mathcal R(A\subset M)=\mathcal R(G \curvearrowright A)$, the right hande side being the p.m.p. equivalence relation generated by the action of $G$ on $A$. Thus, within the class $\mathcal{C}_q$, the objects $M=A\rtimes \Gamma_q(G,K)$ remember the OE class of the action $G \curvearrowright A$, up to stable isomorphism, and hence are partially classified by these OE classes.
\par On the other hand, we can construct a slightly different type of generalized q-gaussians $A \rtimes \Gamma^1_{q}(G,K)=(A\bar{\ten} \Gamma(\ell^2(G)\ten K))\rtimes G$ having the property that $A\rtimes \Gamma^1_{q}(G,K)\cong B\rtimes \Gamma^1_{q}(G',K')$ implies that $G \curvearrowright A$ and $G' \curvearrowright B$ are stably OE and if $G\curvearrowright A$ and $G'\curvearrowright A$ are OE then the associated objects are isomorphic. Hence the classification problem for these objects is almost reduced to the orbit equivalence of the actions (see section 7 for more details). At the time of this writing it is not clear whether orbit equivalence of the action implies isomorphism of the generalized $q$-crossed products in full generality, unless $q=0$. This leaves open the possibility that $A\rtimes \Gamma_q(G,K)$ remembers $q$. The partial converse we can prove is that if $\mathcal R(G \curvearrowright A)\cong \mathcal R(G' \curvearrowright A)$, then $A\rtimes \Gamma^1_{q}(G,K)\cong A\rtimes \Gamma^1_{q}(G',K)$ if the representation is given by conjugation. To be more precise, we have
\begin{theorem}\label{real} Let $A$ be abelian, $|q|<1$ and $K$ infinite dimensional. If $\mathcal R(G\curvearrowright A)=\mathcal R(\tilde{G}\curvearrowright A)$ then $(A \bar{\ten} \Gamma_q(\ell_2(G) \ten K))\rtimes G$ and
$(A \bar{\ten} \Gamma_q(\ell_2(\tilde{G}) \ten K))\rtimes \tilde{G}$ are isomorphic. Conversely, if
 \begin{enumerate}
 \item[i)] $A$ and $\tilde{A}$ are abelian, the inclusions $A \subset (A \bar{\ten} \Gamma_q(\ell_2(G)\ten K))\rtimes G$ and $\tilde{A} \subset (\tilde{A} \bar{\ten} \Gamma_q(\ell_2(\tilde{G}),K))\rtimes G$ are rigid;
 \item[ii)] One of the conditions in Corollary 6.4 holds;
 \item[iii)] $[G,G]$ is ICC and the action of $G$ is free and ergodic,
 \end{enumerate}
then $(A \bar{\ten} \Gamma_q(\ell_2(G)\ten K))\rtimes G \cong (\tilde{A} \bar{\ten} \Gamma_q(\ell_2(\tilde{G})\ten K))\rtimes \tilde{G}$ implies that $\mathcal R(G\curvearrowright A)$ and $\mathcal R(\tilde{G}\curvearrowright \tilde{A})$ are stably isomorphic.
\end{theorem}
Taking $\pi$ to be a \emph{unitary} representation on the \emph{complex} $\ell^2(G)$, we have the following
\begin{cor}Let $A$ be abelian, $|q|<1$ and $K$ infinite dimensional. If $\mathcal R(G \curvearrowright A)=\mathcal R(\tilde{G} \curvearrowright A)$ and $\pi:G \to \mathcal{U}(\ell^2(G))$ is the unitary representation given by conjugation on the complex Hilbert space $\ell^2(G)$ then $A\rtimes \Gamma_q^{\pi}(G,K)$ and $A\rtimes \Gamma_q^{\pi}(\tilde{G},K)$ are isomorphic.
\end{cor}
\par Finally, let's say a couple of words about the proof of Theorem 1.1. The main ideas go back to \cite{PoBe} and, to a lesser extent, \cite{IPP} and \cite{PoI}. The ingredients of the proof are the rigidity of the inclusions $A \subset M=A\rtimes_{\si} \Gamma_q^{\pi}(G,K), B \subset M=B\rtimes_{\rho} \Gamma_q^{\pi}(G',K')$ and the Haagerup property of the groups $G,G'$, together with the existence of two 1-parameter groups of automorphisms of $\tilde{M}_A=A\rtimes_{\si} \Gamma_q^{\pi}(G,K \oplus K)$, $\tilde{M}_B=B\rtimes_{\rho} \Gamma_q^{\pi}(G',K' \oplus K')$, respectively (one for each decomposition), all exploited in a manner which has by now become standard (see e.g. \cite{PoBe, PoI, IPP}). It should be mentioned that in the case of trivial representation $\pi:G \to \mathcal O(\ell_{\rz}^2(G))$ we can handle a much larger class of groups but only by using the recent strong results of Popa and Vaes \cite{PoVaI,PoVaII}.
\par {\bf Step 1.} Let's denote by $\alpha_t^{A}, \alpha_t^{B}$ the one parameter groups of automorphisms associated with the two decompositions. Due to the rigidity of the inclusion $A \subset B\rtimes \Gamma_q(G',K')$, $\alpha_t^{B}$ has to converge uniformly on the unit ball of $A$, which implies that a corner of $A$ embeds into $B\rtimes [G',G']$ inside $M$.
\par {\bf Step 2.} Using the rigidity of the inclusion $A \subset (B\bar{\ten}\Gamma_q(\ell^2_{\rz}(G)\ten K'))\rtimes G'$ together with the Haagerup property, we see that actually a corner of $A$ has to embed into $B$ inside $M$, i.e. $A \prec_M B$.
\par {\bf Step 3.} By symmetry, we also have $B \prec_M A$. Note that we cannot deduce that $A$ and $B$ are unitarily conjugate, as $A,B$ are not MASAs.
\par {\bf Step 4.} Theorem 3.3 in \cite{MeVa} allows us to conclude that $\mathcal R(A \subset M)$ is stably isomorphic to $\mathcal R(B \subset M)$. Since by a separate argument we also have that $\mathcal R (A\subset M)=\mathcal R(G\curvearrowright A)$ and $\mathcal R(B \subset M)=\mathcal R(G'\curvearrowright B)$, we arrive at our conclusion.
\par Throughout the paper we use standard notation in von Neumann algebra theory, see e.g. \cite{Tak}.

{\bf Acknowledgement:} The first author would like to thank Adrian Ioana for many helpful conversations about the results in Section 5, as well as for pointing out his results about rigidity of actions on measure spaces. In addition we thank Stefaan Vaes for fruitful  conversations, and in particular for  bringing \cite{MeVa} to our attention.

\section{Popa's Intertwining Techniques}
We will briefly review the concept of intertwining two subalgebras inside a von Neumann algebra, along with the main technical tools developed by Popa in \cite{PoBe,PoI}. Given $N$ a finite von Neumann algebra, let $P\subset fNf$, $Q\subset N$ be diffuse subalgebras for some projection $f\in N$. We say that \emph{a corner of $P$ can be intertwined into $Q$ inside $N$} if there exist two non-zero projections $p\in P$, $q\in Q$, a non-zero partial isometry $v\in pNq$, and a $*$-homomorphism $\psi:pPp \rightarrow qQq$ such that $v\psi(x)=xv$ for all $x\in pPp$. Throughout this paper we denote by $P \prec_{N} Q$ whenever this property holds, and by
$P \nprec_{N} Q$ its negation. The partial isometry $v$ is called an intertwiner between $P$ and $Q$.

Popa established an efficient criterion for the existence of such intertwiners (Theorem 2.1 in \cite{PoI}). Particularly useful in concrete applications is the following \emph{analytic} description of absence of intertwiners.

\begin{theorem}[Corollary 2.3 in  \cite{PoI}]\label{intertwining-techniques} Let $N$ be a von Neumann algebra and let  $P\subset
fNf$, $Q\subset N$ be diffuse subalgebras for some projection $f \in N$. Then the following are equivalent:
\begin{enumerate}
\item $P \nprec_N Q$.
\item For every finite set $\mathcal F\subset fNf$ and every $\epsilon>0$ there exists a unitary $v\in \mathcal U(P)$ such that
\begin{equation*}\sum_{x,y\in\mathcal F} \|E_Q(xvy^*)\|^2_2 \leq \epsilon. \end{equation*}  \end{enumerate}
\end{theorem}

\begin{defi}Let $(M, \tau)$ be a finite von Neumann algebra, $A \subset M$ a von Neumann subalgebra and $\Phi:M \rightarrow M$ a normal, completely positive, sub-unital, sub-tracial map. We say that $\Phi$ is compact over $A$ if the canonical operator $T_{\Phi}:L^2(M) \rightarrow L^2(M)$ ($T_{\Phi}(\hat{x})=\widehat{\Phi(x)}, x \in M$) belongs to the compact ideal space of $\langle M, e_A \rangle$ (see \cite{PoBe},1.3.3 and \cite{OPCartanI}, 2.7)
\end{defi}

The following result is Prop.2.7 in \cite{OPCartanI}.
\begin{prop}Let $(M,\tau)$ be a finite von Neumann algebra and let $A, P \subset M$ be two von Neumann subalgebras. Let $\Phi:M \rightarrow M$ be a normal, completely positive, sub-unital, sub-tracial map which is compact over $A$ and assume that
\[\inf_{u \in \mathcal{U}(P)} \|\Phi(u)\|_2 > 0. \]
Then $P \prec_M A$.
\end{prop}

\section{Equivalence relations associated to abelian non-maximal abelian subalgebras}
In \cite{MeVa}, Meesschaert and Vaes defined the generalized equivalence relation associated to an inclusion $A=L^{\infty}(X) \subset M$, where $M$ is a type $II_1$ factor and $A$ a diffuse abelian subalgebra of $M$ which is not maximal abelian. This equivalence relation, denoted by $\mathcal R(A \subset M)$, is defined as the measurable equivalence relation on $X$ generated by the graphs of all the partial automorphisms of $X$ associated to the partial isometries $u \in M$ such that $uu^*,u^*u \in A'\cap M, uAu^*=A$. Note that in the case of $A$ being a MASA, this coincides with the standard p.m.p. equivalence relation defined by Feldman and Moore. The following is Theorem 3.3 in \cite{MeVa}:
\begin{theorem}Let $M$ be a type $II_1$ factor with separable predual. Let $A,B \subset M$ be abelian, quasi-regular von Neumann subalgebras satisfying $\mathcal Z(A'\cap M)=A$ and $\mathcal Z(B'\cap M)=B$. If $A \prec_M B$ and $B \prec_M A$, then the equivalence relations $\mathcal R(A \subset M)$ and $\mathcal R (B \subset M)$ are stably isomorphic.
\end{theorem}
\par Just like in \cite{MeVa}, for $A,B$ two abelian von Neumann algebras, PIso$(A,B)$ will denote the set of all partial isomorphisms from $A$ to $B$, that is isomorphisms $\alpha:Aq \to Bp$, where $q \in \mathcal P(A), p \in \mathcal P(B)$. PAut$(A)$ will be used instead of PIso$(A)$. To every $\alpha$ in PIso$(A,B)$ one can associate an $A-B$ bimodule $\mathcal H(\alpha)$ given by $\mathcal H(\alpha)=L^2(Bp)$ and $a\xi b=\alpha(aq)\xi bp$.
\par We will also need the following result, which is Lemma 3.4 in \cite{MeVa}.
\begin{prop}Let $(M,\tau)$ be a tracial von Neumann algebra and $A=L^{\infty}(X,\mu) \subset M$ an abelian von Neumann subalgebra such that $\mathcal Z(A'\cap M)=A$. Let $\mathcal F \subset M$ be a subset such that
\begin{itemize}
\item $M=(\mathcal F \cup \mathcal F^* \cup (A'\cap M))''$;
\item the $\|\cdot\|_2$-closed span of $A\mathcal F A$ is isomorphic, as an $A-A$ bimodule, to a direct sum of bimodules of the form $\mathcal H(\alpha_n)$, where $\rm \alpha_n \in PAut(A)$.
\end{itemize}
Choose nonsingular partial automorphisms $\phi_n$ of $(X,\mu)$ such that $\alpha_n=\alpha_{\phi_n}$ for all n. Then $\mathcal R(A\subset M)$ is generated, up to measure zero, by the graphs of the partial isomorphisms $\phi_n$.
\end{prop}

\section{The Generalized $q$-gaussian algebras}

\subsection{Background on $\Gamma_q(\mathcal H)$}
Let us first recall (see \cite{BoSpe, BoKuSpe}) that for every $-1<q<1$ there is a functor $\Gamma_q$ from the category of real Hilbert spaces with real contractions to the category of finite von Neumann algebras with normal, tracial, ucp maps having the following properties:
\begin{enumerate}
\item For every real Hilbert space $\mathcal H$ there exists a finite von Neumann algebra $\Gamma_q(\mathcal H)$ and a linear map $s_q:\mathcal H \to \Gamma_q(\mathcal H)_{sa}$  such that
 \begin{equation}\label{mom}  \tau(s_q(h_1)\cdots s_q(h_m))\lel \sum_{\si\in P_2(m)} q^{{\rm cr}(\si)} \prod_{\{i,j\}\in \si} (h_j,h_j) \pl
 \end{equation}
and $\Gamma_q(\mathcal H)$ is generated by the $s_q(h)$'s with $h \in \mathcal H$. Here $P_2(m)$ stands for the set of pair partitions of the set $\lbrace 1, \ldots, m \rbrace$ and cr$(\sigma)$ denotes the number of crossings of the pair partition $\sigma$. Sometimes we will drop the subscript $q$ when it's clearly understood from the context and just write $s(h)$ instead of $s_q(h)$.
\item The functor $\Gamma_q$ gives rise to a group homomorphism  $\Gamma_q: \mathcal O(\mathcal H) \rightarrow {\rm Aut}(\Gamma_q(\mathcal H))$ such that
    \[ \Gamma_q(o)(s(h)) = s(o(h)), h \in \mathcal H.\]
\item Let $\mathcal H \subset \mathcal K$ be an inclusion of real Hilbert spaces. Let $P_{\mathcal H}: \mathcal K \rightarrow \mathcal H$ be the orthogonal projection. Then $\Gamma_q(\mathcal H) \subset \Gamma_q(\mathcal K)$ and moreover $E_{\Gamma_q(\mathcal H)}=\Gamma_q(P_{\mathcal H})$, where $E_{\Gamma_q(\mathcal H)}$ denotes the canonical conditional expectation.
\item The von Neumann algebra $\Gamma_q(\mathcal H)$ is represented in standard form on
 \[ L^2(\Gamma_q(\mathcal H))\cong \mathcal{F}_q(\mathcal H)\lel \bigoplus_{n=0}^{\infty} \mathcal{H}_q^{\ten n} \pl ,\]
where $\mathcal{H}_q^{\ten n}$ is the completion of the $n$-fold tensor product of $\mathcal H \ten \cz$ equipped with the inner product
 \[ (h_1\ten \cdots \ten h_n,k_1\ten\cdots \ten k_n)_q
 \lel \sum_{\si\in S_n} q^{{\rm inv(\si)}} (h_{\si(1)}\ten \cdots \ten h_{\si(n)},k_1\ten \cdots \ten k_n) \pl .\]
Here $\mathcal{H}_q^{\ten 0} = \mathbb{C}\Omega$, where $\Omega$ is the vacuum vector. Also the trace on $\Gamma_q(\mathcal H)$ is given by $\tau(x)=(x\Omega,\Omega), x \in \Gamma_q(\mathcal H)$.
For $q=0$ we have the usual inner product in the $(\mathcal H \ten \cz)^{\ten n}$.
\item The formula $\Gamma_q(o)(s(h)) \lel s(o(h))$ can be extended to real contractions $v: \mathcal H \to \mathcal H$. Below we briefly describe how to do this.
 For every real contraction $v$, we have an orthogonal transformation of $\mathcal H \oplus \mathcal H$ given by
 \[    o \lel \kla \begin{array}{cc} v &\sqrt{1-vv^*}\\
                               -\sqrt{1-v^*v} & v*\end{array} \mer      \pl .\]
Then we may define $\Gamma_q(v) \lel E_{\Gamma_q(\mathcal H)} \circ \Gamma_q(o) \circ \iota_{\mathcal H}$, where $\iota_{\mathcal H} : \mathcal H \rightarrow \mathcal H \oplus \mathcal H$, $\iota_{\mathcal H} (h) = (h,0)$. In order to show that $\Gamma_q(v_1 v_2)= \Gamma_q(v_1) \Gamma_q(v_2)$, we have to use the Fock space description. The automorphism $\Gamma_q(o)$ is implemented by $\pi(o) = \oplus_n (o^{\ten n})$ so that
 \[ \Gamma_q(o)(T) \lel \pi(o)T\pi(o^*), T \in \Gamma_q(\mathcal H \oplus \mathcal H) \pl .\]
Similarly the the conditional expectation $E$ commutes with the natural grading. Then the conditional expectation satisfies
 \[ E(\xi) \Om \lel  \oplus_n E^{\ten n}(\xi_n) \]
where $\xi = \oplus_n \xi_n$ is the decomposition in the Fock space. From this it follows that for every contraction
 \[ \Gamma_q(v)(\xi)\Om \lel  (v^{\ten n}\xi_n)_{n\gl 0} \pl .\]
Using this description in $L_2(\Gamma_q(\mathcal H))$ and the injectivity of the inclusion $\Gamma_q(\mathcal H)\subset \mathcal{F}_q(H)\lel L_2(\Gamma_q(\mathcal H))$ it is then easy to deduce that $\Gamma_q$ is a group homomorphism.
\item The most prominent example of such a ucp map arising from a contraction is given by the semigroup of completely positive maps $T_t=\Gamma_q(e^{-t}Id)$. It follows immediately that the generator of $N$ of this semigroup, i.e. $T_t=e^{-tN}$, corresponds to  the usual number operator
 \[ N(\xi_n) \lel n \xi_n, \xi_n \in \mathcal{H}^{\ten n}  \]
on the q-Fock space. Note that a dilation by automorphism $\al_{\theta} \in Aut(\Gamma_q(\mathcal H \oplus \mathcal H))$ is  ``built in'' the construction. Indeed, let $e^{-t}=\cos(\theta)$ and
 \[ o_{\theta} \lel \kla \begin{array}{cc} \cos(\theta) & \sin(\theta)\\
 -\sin(\theta) & \cos(\theta)\end{array}\mer \pl .\]
We denote by $\al_{\theta}=\Gamma_q(o_{\theta})$ and observe that $T_t=E_{\Gamma_q(\mathcal H)} \circ \al_{\theta}|_{\Gamma_q(\mathcal H)}$.
\item For every tensor $\xi \in \mathcal{H}^{\ten n}$ there is an unique element $W(\xi) \in \Gamma_q(\mathcal H)$ (called the \emph{Wick word} for $\xi$) such that $W(\xi)\Omega=\xi$. Due to functoriality, for every real contraction $u: \mathcal H \rightarrow \mathcal H$, we have
\[ \Gamma_q(u)(W(h_1 \ten \ldots \ten h_m)) = W(u(h_1) \ten \ldots \ten u(h_m)), h_1, \ldots, h_m \in \mathcal H \]
\item A concrete description of $\Gamma_q(\mathcal H)$ is given by $\Gamma_q(\mathcal H) = \lbrace s_q(h) : h \in \mathcal H \rbrace '' \subset \mathcal B(\mathcal{F}_q(\mathcal H))$, where for real $h \in \mathcal H$ we have
 \[ s_q(h) \lel l_q(h)+l_q(h)^* \pl ,\]
where $l_q(h)(h_1\ten \cdots h_n)=h\ten h_1\ten \cdots \ten h_n$ is the  creation operator and
 \[ l_q(h)^*(h_1\ten \cdots \ten h_n)
 \lel \sum_{j=1}^m q^{j-1} (h,h_j) h_1\ten \cdots h_{j-1} \ten  \hat{h}_j \ten h_{j+1}\ten \cdots \ten h_n \]
is the adjoint with respect to the $q$-inner product. Here $\hat{h}_j$ means that this vector is omitted. For our analysis it will be important to note that the real linear map $s$ admits a complex extension, also denoted by $s$ to $\mathcal H_{\cz} = \mathcal H\ten_{\rz}\cz\cong \mathcal H \oplus \mathcal H$ (the complexification of $\mathcal H$) given by
 \[ s(h_1+ih_2) \lel s(h_1)+is(h_2) \pl .\]
\end{enumerate}
We will also need an ultraproduct approach to constructing Wick words.
Let us fix $\mathcal{H}$ and $\nen$ and denote by $e_j$ the unit vectors in $\ell_2^n$. By functoriality (1), we see that
 \[ u_n(s_q(h)) \lel \frac{1}{\sqrt{n}}\sum_{j=1}^n s_q(h\ten e_j)  \]
extends to a $^*$-homomorphism from $\Gamma_q(\mathcal{H})$ to
$\Gamma_q(\ell_2^n(\mathcal{H}))$. In particular, we have
 \[ u_n(s_q(h_1)\cdots s_q(h_m))
 \lel \frac{1}{n^{m/2}} \sum_{1 \leq j_1,...,j_m \leq n} s_q(h_1\ten e_{j_1})\cdots s_q(h_m\ten e_{j_m})  \pl .\]
We need to recall some notation. For $1\le j_k\le n$ and a partition $\si$ of $\{1,...,m\}$ we write $\lan j_1,...,j_m\ran=\si$ if
 \[ j_r\lel j_s \quad \Leftrightarrow \exists_{A\in \si} : r,s\in A \pl. \]
In other words indices coincide if they have the same color given by the coloring of $\si$.  We denote by $P_{1,2}(m)$ the set of partitions which only contain singletons and pairs. Let us define
 \[ x_{\si}^n(h_1,...,h_m) \lel \frac{1}{\sqrt{n}^m} \sum_{\lan j_1,...,j_m\ran=\si} s_q(h_1\ten e_{j_1})\cdots s_q(h_m\ten e_{j_m})  \pl .\]
Then we have
 \begin{equation}\label{dcc}
   u_n(s_q(h_1)\cdots s_q(h_m)) \lel \sum_{\si} x_{\si}^n(h_1,...,h_m) \pl .
   \end{equation}
Fix a free ultrafilter $\om$ on the natural numbers. We will make frequent use of the canonical embedding $u_{\om}:\Gamma_q(\mathcal H) \rightarrow \prod_{n,\om} \Gamma_q(\ell^2_n(\mathcal H))$, given by $u_{\om}(x) = (u_n(x))_n$, and we will sometimes identify $\Gamma_q(\mathcal H)$ with its image in the ultraproduct.
\begin{prop}\label{vanish}  For a partition $\si\notin P_{1,2}$ and $p > 2$ we have
 \[  \|x_{\si}^n(h_1,...,h_m)\|_p \kl c(p) n^{1/p-1/2} \pl . \]
\end{prop}

\begin{proof} We need the M\"obius inversion formula for functions $f:\{1,...,n\}\to V$, $V$ a vector space, as it is presented in \cite{Port}. For a partition $\si$ of $\{1,...,m\}$ we define
 \[ \lan \si \ran \lel \sum_{\lan j_1,...,j_d\ran=\si} f(j_1,...,j_d) \] and
  \[ [\si] \lel \sum_{\nu \gl \si} \lan \nu \ran \pl .\]
Let us recall that $\si\le \nu$ if every set (block)  $A\in \si$  is contained in some block of $\nu$. The M\"obius inversion formula states that conversely
 \[ \lan \si \ran \lel \sum_{\nu \gl \si} \mu(\nu,\si) [\nu] \pl \]
holds for some universal (integer valued) function $\mu$. Let us now fix $1\le p\le \infty$ and elements $x_j(1),...,x_j(m)\in L_{mp}(N,\tau)$ and a given partition $\si$. Let $A\in \si$ be a partition with three elements $A=\{k_1,k_2,k_3\}$. For all the other partitions we apply Pisier's unitary trick and find $g_j(k)$ so that
 \[ \tau(g_{j_2}(2)\cdots g_{j_{k_2-1}}(k_2-1)g_{j_{k_2+1}}(k_2+1)
 \cdots g_{j_m}(m)) \lel \begin{cases} 1 & \lan j_1,...,j_m\ran =\si \\
                                       0 & \mbox{else} \end{cases} \pl .\]
This allows us to  define $X_j(k)=g_j(k)\ten x_j(k)$ and write
 \begin{align*}
  [\si] &=  \sum_j x_{j}(1)(\sum_{j_2,...,j_{k_2-1}}X_{j_2}\cdots X_{j_{k_2-1}}(k_2-1))x_j(k_2)\\
  &\quad \quad
 (\sum_{j_{k_2+1},...,j_{k_3-1}}X_{j_{k_2+1}}\cdots X_{j_{k_3-1}}(k_3-1))x_j(k_3)(\sum_{j_{k_3+1},...,j_m}X_{j_{k_3+1}}(k_3+1)\cdots X_{j_m}(m) \\
 &= \sum_j x_j(1)ax_j(k_2)bx_j(k_3)c \\
 &= (\sum_j e_{1,j}\ten x_{j}(1)) (\sum_j e_{jj}\ten ax_j(k_2)b)
 \sum_j e_{j,1}\ten x_j(k_3)c
   \pl .
 \end{align*}
According to H\"{o}lders inequality we find
 \begin{align*}
  \|[\si]\|_p &=  \|(\sum_j x_j(1)x_j(1)^*)^{1/2}\|_{pm}
 (\sum_{j=1}^n \|ax_j(k_2)b\|_q^q)^{1/q} \|(\sum_j x_j(k_3)^*x_j(k_3))^{1/2}\|_{pm}  \|c\|_v \\
 &\kl  \|(\sum_j x_j(1)x_j(1)^*)^{1/2}\|_{pm} \|a\|_{r_1} (\sum_j x_j(k_2)\|_{pm}^q)^{1/q} \|b\|_{r_2} \|(\sum_j x_j(k_3)^*x_j(k_3))^{1/2} \|c\|_v \pl .
 \end{align*}
Here we need $\frac{1}{v}+\frac{1}{r_1}+\frac{1}{r_2}+\frac{2}{pm}=\frac{1}{p}$ and
$\frac{1}{q}=\frac{1}{r_1}+\frac{1}{r_2}+\frac{1}{pm}$. Moreover, according to Pisier's estimate for the $g_j(k)$ we have
 \[ \|a\|_{r_1}\kl C^{r_1} \prod_{1<k<k_2} \tilde{S}(k)  \pl,
 \]
where  $\tilde{S}(k)=\max\{\|\sum_j x_j(k)\|,\|(\sum_j x_j(k)^*x_j(k))^{1/2}\|, \|(\sum_j x_j(k)x_j(k)^*)^{1/2}\|\}$. Thus we
obtain
 \[ \|[\si]\|_p \kl  n^{1/q}
 \sup_j  \|x_j(k_2)\|_{pm} \pl
 C^{m-2} \prod_{k\neq k_2} \tilde{S}(k) n^{1/q}
  \pl .\]
Note here that $q\gl p$. In our situation, we have
 \[ \|\sum_{j=1}^n  s_q(h\ten e_j)\|_{pm} \lel \sqrt{n} \|s_q(h)\|_{pm} \]
by the rotation invariance. For the square function we observe that for $pm\gl 2$ we have
 \[ \|\sum_{j} s_q(h\ten e_j)^2\|_{pm/2}\kl \sum_j \|s_q(h\ten e_j)\|_{pm}^2 \kl n \|s_q(h)\|_{pm}^2 \pl .\]
Thus $\tilde{S}(k)\kl c(pm)\sqrt{n}$ for all $k$. Thanks to our normalization factor we deduce
 \[  \|x_{\si}^n(h_1,...,h_m)\|_p \kl C(pm) n^{1/p-1/2} \pl .\]
For $p>2$ this converges to  $0$. \qd

Note similar arguments can be found in  \cite{Av}. Now we fix a partition $\si=\si_{s}\cup \si_{p}$ of singletons and pairs and $s=|\si_s|,p=|\si_p|$. It follows from \eqref{mom} that
 \begin{align}\label{ipp}
  \lim_{n} \tau(x_{\si'}^n(h'_1,...,h'_m)^*x_{\si}^n(h_1,...,h_m))
  &=   \delta_{s,s'} f_{\si'}(h'_1,...,h'_{m'}) f_{\si}(h_1,...,h_m)\\ \nonumber
  &\quad \quad
 \lim_n  \tau(x_{\si_{s'}}^n(h'_1,...,h'_{m'})^* x_{\si_{s}}^n(h_1,...,h_m))
 \pl,
  \end{align}
where the meaning of the notation will be explained below. Indeed, when calculating the trace, we are only supposed to use pair partitions. If two indices are already combined in  $x_{\si}^n$ we cannot connect them to either a singleton or a pair in $x_{\si'}^n$ because this would produce a partition containing a set with three or four elements, and the result converges to $0$ for $n\to \infty$, according to Proposition 3.1. Therefore we have to connect the singletons of $x_{\si}^n$ with the singletons of $x_{\si'}^n$. Let us put a marker between $x_{\si}^n(h_1,...,h_m)^*$ and $x_{\si'}^n(h_1,...,h_m)$. Each of the singletons from  $x_{\si}^n(h_1,...,h_m)^*$ will be connected to exactly one of the singletons from $x_{\si'}^n(h_1',...,h_m')$ and therefore the left leg of the connection will be left of the marker and in particular left of any pair in $x_{\si}^n$ which crosses over a singleton in $x_{\si}^n$. In other words for
any pair $\{l,r\}\in A\in \si'$ and a singleton $l<i<r$ this will produce exactly one crossing. Therefore we define
 \begin{equation}\label{fff}
  f_{\si}(h_1,...,h_m) \lel q^{{\rm cr}(\si)}\prod_{\{l,r\}\in \si_{p}}  (h_l,h_r) \pl ,
  \end{equation}
where the crossings are counted as usual for pairs and a crossing between a singleton and pair is explained above. Let us not forget that in calculating the inner product we have to sum over all partitions counting the singletons. For this we define
 \[  x_{\si_s}^n(h_1,...,h_m) \lel x_{\{1\},....,\{s\}}(h_{k_1},...,h_{k_s}) \pl ,\]
given by the trivial partition and the indices $\si_s=\{\{k_l\}|l=1,..,s\}$ respecting the order in which they appear in the list of the $h_1,...,h_m$.  The inner product $\tau(x_{\si'_{s}}^n(h'_1,...,h_{m'})^* x_{\si_{s}}^n(h_1,...,h_m))$ then accounts for the missing partitions between singletons. Let us now denote by $X\subset \prod_{n,\om}L_2(\Gamma_q(\ell_2^n(\mathcal{H}))$
the norm-$\|\cdot\|_2$ closed span of the words
  \[ x_{\si}(h_1,...,h_m) \lel (x_{\si}^n(h_1,...,h_m))^{\bullet}, m \geq 0, \si \in P_{1,2}(m), h_1,...,h_m \in \mathcal H \pl .\]

\begin{prop}\label{wick} $X =  u_{\om}(L_2(\Gamma_q(\mathcal{H})))$. More precisely, the decomposition into eigenvectors of the semigroup given by the number operator is given by
 \[ u_{\om}(s_q(h_1)\cdots s_q(h_m)) \lel \sum_{\si\in P_{1,2}(m)}x_{\si}(h_1,...,h_m) \pl.  \]
The Wick words $x_{\si}(h_1,...,h_m)$ belong to $\Gamma_q(\mathcal{H})$ and satisfy $T_t(x_{\si})=e^{-t|\si_s|}x_{\si}$.
\end{prop}

\begin{proof} It follows from \eqref{ipp} that $x_{\si}(h_1,...,h_m)=f_{\si}(h_1,...,h_m)x_{\{1\},...,\{s\}}(h_{i_1},...,h_{i_s})$ where the collection $\{\{i_1\},...,\{i_s\}\}$ is the collection of singletons of $\si$. Since by definition $X$ is the span of the words $x_{\si}$, and by our construction $u_{\om}(s_q(h_1)\cdots s_q(h_m))$ is in $X$, our assertion follows from proving that $x_{\{1\},...,\{s\}}(h_1,...,h_m)$ is in $L_2(\Gamma_q(\mathcal{H})))$ (because $x_{\{1\},...,\{s\}}(h_1,...,h_m)$ is automatically bounded for $|q|<1$). We may prove this by induction on $s$, the number of singletons in $\si$. This is clear for $s=0$ and $s=1$. Now we proceed by induction. Using \eqref{dcc} and \eqref{ipp}, we know that
 \[u_{\om}(s_q(h_1)\cdots s_q(h_m))
 \lel x_{\{1\},....,\{m\}}(h_1,...,h_m)
 + \sum_{|\si_s|<m } f_{\si}(h_1,...,h_m) x_{\si_s}(h_1,...,h_m) \pl .\]
By induction hypothesis the second sum belongs to $u_{\om}(L_2(\Gamma_q(\mathcal H)))$. Hence taking the difference completes the proof. We are left to prove that the $x_{\si}$ are eigenvectors. Indeed, for fixed $t$ we may consider the spectral decomposition of selfadjoint operator $T_t$ on $L_2$ and the projection $q$ onto the orthogonal complement of the eigenspaces for the eigenvalues $\{e^{-tk}:k\in \nz_0\}$. If we were to know that for  $q(h)=h$ of norm $1$, we can approximate $\|h-\sum_j h_j\|_2<1$ so that $T_t(h_j)=e^{-tk_j}h_j$. Then orthogonality implies  $\|h\|^2\le (h-\sum_j h_j,h-\sum_j h_j)<1$, and hence leads to a contradiction.
Thus, for the discrete spectrum of $T_t$ it suffices to show that the $x_{\si}$'s are eigenvectors. Let us  calculate the action of $\al_{\theta}$ on each term $x_{\si}(h_1,...,h_m)$. First we note that the ultraproduct construction also works for $\mathcal{H} \oplus \mathcal{H}$ and commutes with the $(\al_{\theta}^n)^{\bullet}$ action applied component-wise in the ultraproduct. In particular we may apply  $\al_{\theta}^n$ to a word $x_{\si}$
 \begin{align*}
  \al_{\theta}(x_{\si}(h_1,...,h_m) &=
  (\al_{\theta}^n x_{\si}^n(h_1,...,h_m))^{\bullet}
  \lel (x_{\si}^n(o_{\theta}(h_1),...,o_{\theta}(h_m)) \\
  &= f_{\si}(o_{\theta}(h_1),...,o_{\theta}(h_m)) x_{\si_s}(o_{\theta}(h_1),...,o_{\theta}(h_m)) \pl .
  \end{align*}
A glance at \eqref{fff} shows that
 \[ f_{\si}(o_{\theta}(h_1),...,o_{\theta}(h_m))
 \lel f_{\si}(h_1,...,h_m) \]
because $o_{\theta}$ preserves inner products. Moreover, for $x_{\si_s}$, $\al_{\theta}$ only acts on the positions of the single elements. Taking now the orthogonal projection onto $u_{\om}(L_2(\Gamma_q(\mathcal{H}))$, and assuming in addition that $h_1,...,h_m$ are in $\mathcal{H}\oplus 0$, we deduce from the moment formula that
 \[ \tau(x'_{\si_{s}}(k_1,...,k_s)^*x_{\si}(o_{\theta}(h_1),...,o_{\theta}h_m))
   \lel   \cos(\theta)^s  \tau(x'_{\si_{s}}(k_1,...,k_s)^*x_{\si}(h_1,...,h_m))
   \pl .\]
Thus $x_{\si_s}(h_1,...,h_m)$ is an eigenvector of the heat semigroup and we have exactly recovered the Fock space structure.\qd

\begin{rem}\label{wickproduct} The proof reveals a convolution structure  for multiplying Wick words which is independent of $q$.  Indeed, let $\xi$ and $\eta$ be Wick words represented as
 \[  W(\xi)W(\eta)
 \lel
 \big(n^{-m+m'/2}\sum_{\lan j_1,...,j_m\ran=\emptyset} s_{j_1}(h_1)\cdots s_{j_m} (h_m)\sum_{\lan j'_1,...,j'_{m'}\ran=\emptyset} s_{j'_1}(h'_1)\cdots s_{j'_{m'}}(h'_{m'})\big )^{\bullet} \pl .\]
Then we may apply Proposition \ref{vanish}, and  hence it remains to sum over singleton/pair partition on $\{1,...,m+m'\}$. Note however, that thanks to the additional condition  $\lan j_1,...,j_m\ran=\emptyset$ one cannot pair singletons from $\xi$ or $\eta$. We end up with the sum of bipartite partitions  for the possible pairs. However, according to the argument before \ref{fff}, the $x_{\si}$ corresponding to every  bipartite  partition $\si\in P_{1,2}(m+m')$  can be replaced by a new Wick word whose length is given by the singletons in $\si$ and the additional coefficient $f_{\si}(\xi, \eta)=q^{{\rm cr(\si)}}\prod (h_i,h'_{i'})$ depending on the pairs. This means
\begin{align} \label{wp}
 W(\xi)W(\eta)
  &= \sum_{\si\in P_{1,2}(m+m')} f_{\si}(\xi,\eta)
   W(\xi_{A(\si)}\ten \eta_{B(\si)}) \pl,
   \end{align}
where $A(\si)$, $B(\si)$ are the unions of singletons of $\si$ corresponding to vectors in $\xi$ and $\eta$, respectively and $\xi_A$, $\eta_B$ are obtained  by erasing the tensors in $\xi$ and $\eta$ whose indices do not belong to $A$ and $B$, respectively. We will freely use the structure of this multiplication, also for products of three elements.
\end{rem}

\subsection{$q$-Gaussian Group Measure Space Construction}

Let's first recall the construction of the Gaussian action associated to a representation (see for example \cite{PeSi}). Let $\pi:G \rightarrow \mathcal{O}(\mathcal H)$ be an orthogonal representation on a real Hilbert space, and consider the abelian von Neumann algebra $(D, \tau) \cong \Gamma_1(\mathcal H)$ generated by a family of unitaries $\omega(\xi),\xi \in \mathcal H$, subject to the following relations:

\begin{enumerate}
\item $\omega(\xi_1)\omega(\xi_2)=\omega(\xi_1+\xi_2)$, for any $\xi_1, \xi_2 \in \mathcal H$;
\item $\omega(-\xi)=\omega(\xi)^*$, for any $\xi \in \mathcal H$;
\item $\tau(\omega(\xi))=\exp(-||\xi||^2)$, for any $\xi \in \mathcal H$.
\end{enumerate}

The \emph{Gaussian action} of $G$ on $(D, \tau)$ is defined by $\tilde{\pi}_g(\omega(\xi))=\omega(\pi_g(\xi))$, for all $g \in G$ and $\xi \in \mathcal H$.
This construction can be generalized for $q \neq 1$. Namely, let $\pi : G \rightarrow \mathcal O (\mathcal H)$ be an orthogonal representation of $G$ on the real Hilbert space $\mathcal H$. Then the $q$-gaussian action of $G$ on $\Gamma_q(\mathcal H)$ can be defined using the above functoriality properties by
\[ \tilde{\pi}_g(s(h)) = s(\pi_g(h)), \forall g \in G, h \in \mathcal H. \]

\begin{rem} \label{diag} Let $G$ be a group acting trace preservingly on $A$ via $\si$, and  $\pi:G\to \mathcal O(\mathcal H)$ be an orthogonal representation of $G$ on a real Hilbert space. Then $G$ acts diagonally on $A \bar{\ten} \Gamma_q(\mathcal H)$ by
 \[ \rho_g (a\ten s(h)) \lel \si_g(a) \ten s(\pi_g(h)) \]
\end{rem}

\begin{lemma} The semigroup of completely positive maps $T_t$ given by the number operator on $\Gamma_q(H)$ extends to the crossed product  $(A \bar{\ten} \Gamma_q(\mathcal H)) \rtimes_{\rho} G$ and admits a factorization
 \begin{equation}\label{dialat}
 T_t \lel E_{(A \bar{\ten} \Gamma_q((\mathcal H,0))\rtimes G} \circ (\al_{o_t}\rtimes 1_G)|_{(A \bar{\ten} \Gamma_q((\mathcal H,0))\rtimes G} \pl.
 \end{equation}
The eigenspaces are of the form
  \[ W_n \lel {\rm span}\{(a\ten W(h_1\ten \cdots \ten h_n))u_g: a\in A, h_1\ten \cdots h_n\in H_q^n, g \in G\} \pl ,\]
where the $u_g$'s are the canonical unitaries.
\end{lemma}

\begin{proof} Let us recall that the semigroup $T_t=e^{-tN}$ given by the number operator $N$ is implemented via
 \[ T_t \lel E_{\Gamma_q((\mathcal H,0))} \al_{o_t}|_{\Gamma_q((\mathcal H,0))} \pl ,\]
where  $\Gamma_q((\mathcal H,0))\subset \Gamma_q(\mathcal H \oplus \mathcal H)$ is viewed as a von Neumann  subalgebra, and the orthogonal matrix $o_t=\kla \begin{array}{cc} \cos(\theta)& -\sin(\theta)\\
                  \sin(\theta)& \cos(\theta)\end{array}\mer$
is the rotation with $\cos(\theta)=e^{-t}$. Since $T_t$ is completely positive we may consider $\hat{T}_t=id_A\ten T_t\ten id_{B(\ell_2(G))}$. Moreover, $T_t$ commutes with
the group action $\al_{\pi(g)}$, and this implies that
 \[  \hat{T}_t((A\ten \Gamma_g(\mathcal H))\rtimes G) \subset
 \hat{T}_t((A\ten \Gamma_g(\mathcal H))\rtimes G) \pl .\]
Thus the restriction $T_t\rtimes 1_{G}$ to $(A\ten \Gamma_g(\mathcal H)) \rtimes G$ is well-defined. Similarly, we see that $\al_{o_t}\rtimes 1_G$ is a well-defined automorphism of $\Gamma_q(\mathcal H \oplus \mathcal H)$ using the diagonal action (compatible with second quantization of $\kla \begin{array}{cc} o_t& 0\\
                  0& o_t\end{array}\mer$).
Thus  by restriction, we obtain \eqref{dialat}. Since $T_t \rtimes 1_G $ is trace preserving, we may consider $T_t$ as acting on
 \begin{align*}
 L_2( (A\rtimes \Gamma_q(\mathcal H))\rtimes G) &=
  L_2(A)\ten L_2(\Gamma_q(\mathcal H))\ten \ell_2(G) \\
 &=  L_2(A) \ten \F_q(\mathcal H) \ten \ell_2(G) \\
 &= \sum_{n=0}^{\infty} L_2(A)\ten \mathcal H_{\cz}^{\ten n}\ten \ell_2(G) \pl .
 \end{align*}
Since $T_t\rtimes 1_{G}$ commutes with the group action, we see that on $L_2$ this is just $id\ten T_t\ten id$. Moreover, the eigenspaces for $T_t$ are exactly those spanned by  tensors of a fixed length, i.e. $T_t(a\ten (h_1\ten \cdots \ten h_n)\ten b)=e^{-tn} a\ten (h_1\ten \cdots \ten h_n)\ten b$. It the follows immediately that for $a\in A$ and $g\in G$, the element $w_n(a,h_1,...,h_n,g)=
(a\ten W(h_1\ten \cdots \ten h_n))u_g$ is in
$(A\rtimes \Gamma_q(\mathcal H))\rtimes G$. Clearly, the linear span of elements  $w_n \Om$ is dense in $L_1(A)\ten H_{\cz}^{\ten n}\ten \ell_2(G)$. \qd

\subsection{$q$-gaussian with group action}

We may apply the $q$-functor in particular to the real Hilbert space $\ell_{\rz}^2(G)$ for $G$ a countable discrete group. Consider the orthogonal representation of $G$ on $\ell_{\rz}^2(G)$ given by
\[ \pi_g (\delta_h) = \delta_{ghg^{-1}} \]

\begin{rem}i) We have
  \[ \pi_g(s_q(\delta_h)) \lel s_q(\delta_{ghg^{-1}}) \pl .\]
iii) The  element
 \[ S_q(g) \lel s_q(e_g)\la(g) \]
satisfies $S_q(g)^*\lel S_q(g^{-1})$.
\end{rem}

\begin{proof} i) trivial by second quantization.
For the proof of ii) we write $R=\frac{\la(g)+\la(g^{-1})}{2}$, $I=\frac{\la(g)-\la(g^{-1})}{2i}$ and get
 \[\delta_g \lel \la(g) \delta_e \lel R+iI \pl .\]
Thus
 \[ s_q(\delta_g) \lel s_q(R)+is_q(I)\]
Hence
 \begin{align*}
 (s_q(\delta_g)\la(g))^* &= \la(g^{-1}) (s_q(R)-is_q(I))
 \lel \pi_{g}^{-1}(s_q(R)-i\pi_{g^{-1}}(I))\la(g^{-1}) \\
 &= (s_q(R)-is_q(I))\la(g^{-1}) \pl .
 \end{align*}
Finally we note that
 \[ R-iI\lel \frac{\la(g)+\la(g^{-1})}{2}-\frac{\la(g)-\la(g^{-1})}{2i}
 \lel \la(g^{-1}) \pl .\]
Thus $(R-iI) \delta_e = \delta_{g^{-1}}$ as asserted. \qd

\begin{defi} Let $G$ be a discrete, countable, infinite group and $K$ a separable real Hilbert space. Let $\sigma : G \curvearrowright A$ a trace preserving action on a finite von Neumann algebra $A$, and $\pi : G \rightarrow \ell_{\rz}^2(G)$ an orthogonal representation of $G$. For a subset $F \subset G$, we define
 \[ A\rtimes \Gamma_q^0(F,K) \]
as the von Neumann subalgebra of $(A \bar{\ten} \Gamma_q(\ell^2(G) \ten K)) \rtimes_{\rho} G$ generated by $A$ and $S_q(g \ten k) = s_q(\delta_g \ten k)u_g$, $g \in F$, $k \in K$, where $\rho$ is the diagonal action of $G$ on $A \bar{\ten} \Gamma_q(\ell_{\rz}^2(G) \ten K)$ associated to the action $\sigma : G \curvearrowright A$ and to the orthogonal representation
\[ \pi \ten id : G \rightarrow O(\ell_{\rz}^2(G) \ten K) \pl. \]
For $K=\rz$ we simply write $A\rtimes \Gamma_q^0(F)$.
\end{defi}

\begin{lemma} \label{inf} Let $K$ be infinite dimensional.
\begin{enumerate}
\item[o)] Let $g\in G$ and $k\in K$. Then
  \[ S_q(g \ten k)a \lel \si_g(a)S_q(g \ten k) \]
\item[i)] $T_t$ leaves $A \rtimes \Gamma_q^0(G,K)$ invariant.
\item[ii)] $L^2(A \rtimes \Gamma_q^0(G,K))=\bigoplus_{n \geq 0} X_n$, where $X_n$ is the $||\cdot||_2$-closed linear span of the elements of the form
     \[ x_n \lel aW((\delta_{g_1} \ten k_1) \ten \cdots \ten (\delta_{g_n} \ten k_n))u_g \]
with $a\in A$ and $(g_1\cdots g_n)g^{-1}\in [G,G]$ for $q\neq 0$; moreover $T_t(x_n)=e^{-tn}x_n$. For $q=0$ and $n=1$ we will have $g=1$ and for other values $W(\delta_{\tilde{g}_1} \ten k_1 \ten \cdots \ten \delta_{\tilde{g}_m} \ten k_m)u_{g_1 \cdots g_m}$ (see proof below).
\end{enumerate}
\end{lemma}

\begin{proof} The first relation is easy to check. Indeed, we have
 \begin{align*}
 S_q(g \ten k)a &=  s_q(\delta_g \ten k)u_g a
 \lel (s_q(\delta_g \ten k) \ten 1)(1\ten  \si_g(a))u_g = \si_g(a)S_q(g \ten k)\pl .
 \end{align*}
For the proof of ii) let us start with
 \begin{align*}
 S_q(g_1\ten k_1)\cdots S_q(g_m\ten k_m)
 &= s_q(\delta_{g_1}\ten k_1)u_{g_1}\cdots s_q(\delta_{g_m}\ten k_m)u_{g_m} \\
 &\lel s_q(\delta_{\tilde{g}_1}\ten k_1)\cdots s_q(\delta_{\tilde{g}_m}\ten k_m) u_{g_1}\cdots u_{g_m} \pl ,
 \end{align*}
Here we have $\tilde{g_j}=g_j$ in case of the trivial action, and $\tilde{g_j}=(g_1\cdots g_{j-1})g_j(g_1\cdots g_{j-1})^{-1}$ in case of the conjugation action.  Now we recall our ultra-product procedure to
`extract' the Wick-words from the first term. Indeed, following Proposition \ref{wick} we have
 \begin{align*}
 & s_q(\delta_{\tilde{g}_1}\ten k_1)\cdots s_q(\delta_{\tilde{g}_m}\ten k_m)
  \lel   \sum_{s=1}^m
 \sum_{\si\in P_{1,2},|\si_s|=\{i_1,...,i_s\}}
 f_{\si}(\delta_{\tilde{g}_1}\ten k_1,\cdots,\delta_{\tilde{g}_m}\ten k_m)\\
 &\quad \quad  \quad \quad \quad \quad \quad \quad \quad \quad \quad \quad \quad \quad \quad \quad \quad \quad W((\delta_{\tilde{g}_{i_1}}\ten h_{i_1})\ten \cdots \ten
  (\delta_{\tilde{g}_{i_s}}\ten h_{i_s}))  u_{g_1}\cdots u_{g_m}  \pl .
  \end{align*}
We recall  the words of length $s$ are obtained by choosing $s$ singleton sets and evaluating the factor
 \[  f_{\si}(\delta_{\tilde{g}_1}\ten k_1,\cdots,\delta_{\tilde{g}_m}\ten k_m)
 \lel q^{{\rm cr}(\si)} \prod_{\{l,r\}\in \si} \tau(u_{\tilde{g_l}}u_{\tilde{g_r}}) (h_l,h_r) \]
obtained from `eliminating the pairings in $\si$' and in particular
 \[  f_{\si}(\delta_{\tilde{g}_1}\ten k_1,\cdots,\delta_{\tilde{g}_m}\ten k_m)
 \lel q^{{\rm cr}(\si)} \prod_{\{l,r\}\in \si} \tau(u_{\tilde{g_l}}u_{\tilde{g_r}}) (h_l,h_r) \]
only depends on $\si$ and the vectors appearing in the pair partitions of $\si$. Let us note that we only obtain a non-trivial term $f_{\si}$ if $\tilde{g_l}\tilde{g_r}=1$ for all the pairs $\{l,r\} \in \si$. Let $p:G \to G/[G,G]$ be the canonical homomorphism. Then we deduce from $f_{\si}\neq 0$ that $p(\tilde{g}_l)p(\tilde{g}_r)=1$ and hence $p(\tilde{g}_1 \cdots \tilde{g}_m)=p(g_{i_1}\cdots g_{i_s})$. This shows that $(g_{i_1}\cdots g_{i_s})^{-1}g_1\cdots g_m\in [G,G]$ in both cases. Now we have to show that for $q\neq 0$ all the expressions $W((\delta_{g_1}\ten h_1)\ten \cdots \ten (\delta_{g_n}\ten h_n))u_g$ with $g_1\cdots g_ng^{-1}\in [G,G]$ will occur. We will prove for $n=0$ first. This will be done by considering the projection of the eigenspace $P_0$ of Number operator. In case of the trivial action we note that for orthogonal unit vectors $h_1$, $h_2$ we have
 \begin{align*}
 & P_0(S_q(g_1\ten h_1)S_{q}(g_2\ten h_2)S_q(g_1^{-1}\ten h_1)S_q(g_2^{-1}\ten h_2)))\\
 &=
   P_0(s_q(g_1\ten h_1)s_q(g_2\ten h_2)s_q(g_1^{-1}\ten h_1)s_q(g_2^{-1}\ten h_2)) u_{g_1}u_{g_2}u_{g_1^{-1}}u_{g_2^{-1}} \\
 &=
 q \tau(u_{g_1}u_{g_1}^{-1}) \tau(u_{g_2}u_{g_2}^{-1}) u_{g_1g_2g_1^{-1}g_2^{-1}} \lel q u_{g_1g_2g_1^{-1}g_2^{-1}} \pl .
 \end{align*}
In case of the conjugation action we define  $G_2=
g_1g_2^{-1}g_1^{-1}$ and $G_1=G_2g_1G_2^{-1}$. For  perpendicular unit vectors we deduce from the  $q$-relations that
 \begin{align*}
  &P_0(S_g(\delta_{G_1}\ten h_1)S_q(\delta_{G_2}\ten h_2)S_q(\delta_{g_1}\ten h_1)S_q(\delta_{g_2}\ten h_2)) \\
  &= q \pl \tau(G_1(G_1G_2)g_1(G_1G_2)^{-1})\tau(G_1G_2G_1^{-1}(G_1G_2g_1)g_2(G_1G_2g_3)^{-1}) u_{G_1}u_{G_2}u_{g_1}u_{g_2} \\
  &= q u_{g_1g_2^{-1}g_1^{-1}g_2} \pl .
  \end{align*}
For arbitrary $s$ we choose a sequence of unit vectors $h_n$ and $k_n$ which converge to $0$ weakly. Using the singleton-pair from Proposition \ref{wick} and the weak convergence we deduce that
 \begin{align*}
  &q W((\delta_{g_1} \ten h_1)\ten \cdots \ten (\delta_{g_s}\ten h_s))u_{g_1\cdots g_s}u_v
  \lel \lim_{n\to \infty}  \\
 &\quad
P_s\big(S_q(g_1\ten h_1)\cdots S_q(g_s\ten h_s)S_q(g_{s+1}\ten h_n)
 S_q(g_{s+2}\ten k_n)S_q(g_{s+3}\ten h_n)S_q(g_{s+2}\ten k_n)\big)
 \end{align*}
for a suitable choice of  $g_{s+1},g_{s+2},g_{s+3},g_{s+4}$. This shows that we have $L^2(A\rtimes \Gamma_q^0(G,K))\subset \oplus_{n\gl 0} X_n$ for $K$ infinite dimensional. Now we have to show i). Since $T_t(x_{\si})=e^{-ts}x_{\si}$ it certainly suffices to
show that all words $x_{\si}$ are actually in $L^2(\Gamma_q^0(G,K))$. Let us embed $K\subset K\ten \ell_2$ via $u(k)=k\ten e_0$. For a given partition $\si$ with a pair $B \in \si$, we may introduce orthogonal vector $e_B$ for any pair $B\in \si$. Then we define $e_k=e_B$ if $k\in B$ and $e_k=e_0$ for singletons. Let $u(s(h))=s(h\ten e_0)$ be the canonical embedding. Then we obtain
 \[ u(x_{\si}) \lel P_kE_{\Gamma_q(K\ten e_0)}(s(h_1\ten e_1)\cdots s(h_m\ten e_m)) \pl .\]
We may  replace $P_kE_{\Gamma_q(K\ten e_0)}$ by a suitable limit.
Indeed, for pairs $\{l,r\}=B\in \si$ we see that $h_l\ten e_B$ and $h_r\ten e_B$ are orthogonal to the other variables. Hence  we can find orthogonal transformations $o_n$ such that $o_n(h_l\ten e_B)$ converges to $0$ weakly and $o(h\ten e_0)=(h\ten e_0)$. Using the moment formula we deduce that
 \[ u(x_\si) \lel w^*-\lim \pl  \Gamma_q(o_n)(s(h_1\ten e_1)\cdots s(h_m\ten e_n)) \pl .\]
For finite dimensional $K$ we deduce that Wick words for $K\ten e_1$ are contained in $\Gamma_q^0(G,K\ten \ell_2)$.
For infinite dimensional $K$ we know that $K$ and $K\ten \ell_2$ are isomorphic. Following the procedure above we can find a sequence of vectors $h_j(n)$ which converge to $0$ weakly and such that $h_j(n)=h_{j'}(n)$ holds for $\{j,j'\}=B\in \si$. Then the moment formula shows that
 \[ x_{\si}(S_q(g_1\ten h_1)\cdots S_q(g_m\ten k_m))
  \lel w^*-\lim_n \pl S_q(g_1\ten \tilde{h}_1(n))\cdots S_q(g_m\ten \tilde{h}_m(n)) \]
where $\tilde{h}_j(n)=h_j$ for a $j$ singleton is a constant,   and $\tilde{h}_j(n)=h_j(n)$ are chosen as above
for the pairs in $\si$. Hence $x_{\si}\in L_2(A\rtimes \Gamma_q^0(G,K))$ for any $\si$ if $K$ is infinite dimensional, and in particular $T_t$ leaves $L_2(A\rtimes \Gamma_q^0(G,K))$ because it is the span of eigenvectors.
\qd

\begin{cor}
\begin{enumerate}
\item[i)]   Let $M=A\rtimes \Gamma_q(F,K)\subset A\rtimes\Gamma_q^0(F,K \ten \ell_2)$ be the von Neumann subalgebra which is invariant under all $\Gamma_q(o)$ with $o|_{K\ten \cz e_1}=id_{K \ten \cz e_1}$. Then $M=A\times \Gamma_q(F,K)$ is invariant under the maps $T_t$ for all $t$, and $L_2(M)$ is the direct sum of  eigenspaces of the number operator.

\item[ii)] Let $C^0\subset M$ be the sub-algebra generated by elements of the form $S_q(g_1)\cdots S_q(g_m)$ with $\prod_j g_j=1$ and $C$ be the von Neumann subalgebra generated by $\bigcup_t T_t(C_0)$. Then $C^0 \subset A' \cap M$ and $C \subset A'\cap M$.
\end{enumerate}
\end{cor}

\begin{proof} We have discussed i) in the proof of Lemma \ref{inf}.
For the proof of ii) we note that $\prod g_j=1$ implies that
 \begin{align*}
 &S_q(g_1\ten k_1)\cdots S_q(g_m\ten k_m)a
  \lel S_q(g_1\ten k_1)\cdots S_q(g_{m-1}\ten h_m) \si_{g_{m-1}}(a)S_q(g_m\ten k_m) \\
 &\lel
  \si_{g_1\cdots g_m}(a) S_q(g_1\ten k_1)\cdots S_q(g_m\ten k_m)
  \lel a S_q(g_1\ten k_1)\cdots S_q(g_m\ten k_m)
  \pl .
  \end{align*}
Thus we deduce the second assertion for $C^0=C^0(K)$. We may apply  also this observation for infinite dimensional $K$. Taking a glance at the  proof of Lemma \ref{inf} we deduce that for infinite dimensional $K$ we still have $T_t(C^0(K))\subset C^0(K)$ because we are not changing the $g_i$'s only the Hilbert space vectors.
Hence $T_t(C^0(K))\subset \overline{C^0(K)}\subset A'$. For finite dimensional $K$, we observe that $C(K)\subset C(K\ten \ell_2)$ and this completes the proof. \qd

For our deformation arguments we need additional estimates for  products of Wick words. In $A\rtimes \Gamma_q(G,K \oplus K)$ we use the notation $E_{K \oplus 0}=E_M$, $E_{0\oplus H}$ for the normal conditional expectations on $A\rtimes \Gamma_q(G,K \oplus 0)=M$, $A \rtimes \Gamma_q(G, 0 \oplus K) = \al_{\pi/2}(M)$, respectively.

\begin{prop}\label{bcontrol} Let us assume that $G$ acts trivially on $\ell_2(G)$. Let $x_1,x_2\in A\rtimes \Gamma_q(G,K \oplus K)$ be Wick words. Then
 \[ V_{x_1,x_2} \lel \{E_{K \oplus 0}(x_1xx_2): x\in A\rtimes \Gamma_q(G,0 \oplus K)\} \]
is contained in a finite dimensional module over $A\rtimes [G,G]$. If $G$ acts by conjugation, $Q\subset A\rtimes \Gamma_q(G,0 \oplus K)$ is rigid and $G$ has  the Haagerup property, then
 \[ V_{x_1,x_2}((Q)_1) \lel \{E_{H\oplus 0}(x_1xx_2): x\in (Q)_1 \} \]
is contained in the $L_2$ closure of a finite module over $A\rtimes[G,G]$, i.e. for every $\eps>$ there exists a finite dimensional right module $H$ over $A\rtimes [G,G]$ such that $V_{x_1,x_2}((Q)_1)$ is contained in an $\eps$ neighborhood of $H$.
\end{prop}

\begin{proof} Let $x_1=a_1W(\xi_1)u_{g_1}$ with $g_1 \in [G,G]$, $x_2=a_2W(\xi_2)u_{g_2}$ and finally $x=aW(\eta)u_h$ with $h\in [G,G]$. Then we find
 \[u_ha_2W(\xi_2) \lel W(\xi_2) a_2u_h \]
and hence
 \begin{align*}
E_{K \oplus 0}(x_1xx_2)&= a_1 E_{K \oplus 0}(W(\xi_1)W(\eta)W(\xi_2)) \si_{g_1}(aa_2)u_{g_1gg_2}
 \end{align*}
Now we apply the procedure of Proposition \ref{wick}, more precisely \ref{wp}. Thus we may rewrite the product of three Wick  words as a linear combination of Wick words. In this process of ``reduction'' we have to combine certain singletons from $\xi$, $\eta$ and $\xi_2$ using pair partitions. Note that thanks to \ref{fff} we known that orthogonal vectors from $\{0\} \oplus K$ and $K \oplus \{0\}$ will be combined through the inner product. Adding the conditional expectation onto $K \oplus 0$ means that all components of $W(\eta)$ have to be paired with components of $\xi_1$ or $\xi_2$ which are in $\{0\} \oplus K$. In particular $\eta$ cannot have a length which exceeds the number of elements from $0 \oplus K$ in $\xi_1$ and $\xi_2$ combined. More precisely, if $\xi_1=h_1\ten \cdots \ten h_m$ and $\xi_2=h'_{1}\ten \cdots \ten h'_{m'}$ we
end up with a finite linear combination of the form
$W(\xi^1_A\ten \xi^2_B)=W(h_{i_1}\ten \cdots h_{i_k}\ten h'_{i_{k+1}}\ten \cdots h'_{i_l})$ where the $i_j$ range through all indices $j$ with $h_j\in K \oplus \{0\}$ or $h'_{j'}\in K \oplus \{0\}$, and hence only depend on our given variables $\xi_1$ and $\xi_2$. The scalar coefficients in front of $W(\xi^1_A\ten \xi^2_B)$ depend on the middle term $W(\eta)$. Thus we find a finite dimensional vector space over $A\rtimes [G,G]$. In the case of the conjugation action, we have to modify the formula and find
 \[ E_{K \oplus 0}(x_1xx_2) \lel a_1E_{K \oplus 0}(W(\xi_1)W(\si_{g_1}(\eta))W(\si_{g_1g_2}(\xi_2)) \si_{g_1g_2}(aa_2)u_{g_1gg_2} \pl .\]
Thus  for the whole algebra $Q=A\rtimes\Gamma_q(G,0 \oplus K)$ we find the linear combinations of all Wick words with combined coefficients from $\xi_1$ and $u_g(\xi_2)$.  Now we assume that \[ Q\subset A\rtimes\Gamma_q(G,0 \oplus K) \subset (A\bar{\ten} \Gamma_q(\ell_2(G)\ten (0\oplus K)))\rtimes G \]
is rigid and $G$ has the Haagerup approximation property. Given $\eps>0$ we can find a finite subset $F\subset G$ such that
  \[ \|x-id\ten P_F(x)\|_2 \kl \eps \|x\| \pl \]
for all $x\in Q$. Thus we may approximate $x$ with linear combination of   Wick words $aW(\xi)u_g$ with $g\in F$. Then we have to consider
all $F$ conjugates of $\si_{g_1g_2}(\xi_2)$ and, as in the first part we obtain a finite dimensional space spanned by  $\xi^1_A\ten \xi^2_B$ coming from $\xi^1$ and $\si_{gg_1g_2}(\xi_2)$ with $gg_1g_2$ running through a finite set. Then we may define the right $A\rtimes [G,G]$ module $H$ spanned by vectors of the form $W(\xi^1_A\ten \si_{gg_1g_2}(\xi_2)_B)u_{g(\xi^1_A\ten \si_{gg_1g_2}(\xi_2)_B)}$, where $g\in F$, $g_1$, $g_2$ fixed. Indeed, our conditional expectation $E_{K \oplus 0}$ does not leave the algebra $A\rtimes \Gamma_q(G,K)$ and for every tensor $\xi=(\delta_{g_1}\ten h_1)\ten \cdots (\delta_{g_n}\ten h_n)$ we have an assigned group elements $g_\xi=g_1\cdots g_n$ so that
the conditional expectation yields elements in the span of $W(\xi)u_{g(\xi)}u_h$ with $h\in [G,G]$. The group elements $u_h$ will depend on the middle expression $W(\eta)u_{g(\eta)h'}$ with $g(\eta)h'\in F$. Let us denote by $P_{H}$ the orthogonal projection onto the right $A\rtimes [G,G]$ module. Then we deduce that for $\|x\|\le 1$ we have
 \begin{align*}
 &\|E_{K \oplus 0}(W(\xi_1)u_{g_1}xW(\xi_2)u_{g_2})-P_H(E_{K \oplus 0}(W(\xi_1)u_{g_1}xW(\xi_2)u_{g_2})\|_2\\
 &\le \|W(\xi_1)u_{g_1}x-(id\rtimes P_F)(x)W(\xi_2)u_{g_2})\|\\
 &\quad
 + \|E_{K \oplus 0}(W(\xi_1)u_{g_1}P_F(x)W(\xi_2)u_{g_2})-P_H(E_{K \oplus 0}(W(\xi_1)u_{g_1}P_F(x)W(\xi_2)u_{g_2}))\|_2 \\
 &\le \eps \|W(\xi_1)\| \|W(\xi)\|  \pl .
 \end{align*}
This concludes the proof in the conjugation case. \qd

\subsection{Factoriality}

We want to show that, under certain conditions, the algebras $C^0$, $C$ and $M = A \rtimes \Gamma_q(F,K)$ defined above are factors.

\begin{prop} \label{Kro}Let $q\gl 0$.  There exists a $k=k(q)$ such that for $|F|\gl k(q)$, the von Neumann algebra $C^0$ is a factor.
\end{prop}

\begin{proof} This follows closely the argument from \cite{Kro}. Indeed, the aim is to show that for certain vectors $x_g\in C^0$, the operator
 \[ T \lel  \sum_{g\in F} |L(x_g)-R(x_g^*)|^2
 \lel \sum_g L(x_g^*x_g)+R(x_g^*x_g)-L(x_g^*)R(x_g^*)-R(x_g)L(x_g) \]
is strictly positive on $L_2(C)\ominus \cz id$. Note that in our case $C\subset \Gamma_q(F,K)$ which allows us to obtain some norm estimates. More specifically, we first consider  $g\in F$ and write $e_g=\al+i\beta$ for orthogonal vectors $\al$ and $\beta$ of norm $1/\sqrt{2}$. Then we see that
 \[ s(g)^*s(g) \lel (s(\al)-is(\beta))(s(\al)+is(\beta))\lel s(\al)^2+s(\beta)^2+ i (s(\al)s(\beta)-s(\beta)s(\al)) \]
and
 \[ s(g)s(g)^* \lel (s(\al)+is(\beta))(s(\al)-is(\beta))\lel s(\al)^2+s(\beta)^2- i (s(\al)s(\beta)-s(\beta)s(\al)) \]
Thus we see that $x_g=s(\al)^2+s(\beta)^2- cid$ belongs to $B$. In the following we will abuse notation and assume that $a=a_g=\sqrt{2}\al$, and $b=b_g=\sqrt{2}\beta$ are unit vectors and omit the factor $2$. Let us note, however, that for $g\neq g'$ the $a_g,b_g$ are both perpendicular  to $a_{g'}$, $b_{g'}$. For normalized $\tilde{a}$ we have $\tau(s(\tilde{a})^2)=1$ and hence will choose $c=2$ above. Then we note that
 \begin{align*}
  &s(a)s(a)(x_g)\Om \lel  s(a)^2[a\ten a+b\ten b] \\
  &=s(a)[a\ten a\ten a+(1+q)a+a\ten b\ten b] \\
  &=(a\ten a\ten a\ten a+ (1+q+q^2)a\ten a +(1+q)a\ten a + (1+q)\Om + a\ten a\ten b\ten b+ b\ten b) \pl .
  \end{align*}
This gives
 \[ (s(a)^2-1)x_g \lel a^{\ten_4}+ a\ten a \ten b\ten b + (1+2q+q^2)(a\ten a)+(1+q)\Om \]
and by symmetry
 \begin{align*}
 x_g^2&= W(a_g^{\ten_4})+W(b_g^{\ten_4}) + W(a_g\ten a_g\ten b_g\ten b_g)+  W(b_g\ten b_g\ten a_g\ten a_g)\\
 &\quad + (1+2q+q^2)(W(a_q\ten a_q)+W(b_g\ten b_g)) +2(1+q)1 \pl.
 \end{align*}
Now, we may apply the Wick word formula from \cite{Nou} and deduce by orthogonality  that
 \[ \|\sum_{g\in F} L(x_g^2)+R(x_g^2)-4(1+q)1\|\kl C_q \sqrt{|F|} \pl .\]
The more challenging part is to understand $L(x_g)R(x_g)$ on $\cz^{\perp}$. Indeed, we have
 \begin{align*}
  R(s(a)^2)\xi &= R(s(a))(\xi\ten a+ r_-(a)\xi)\\
 &= \xi\ten a\ten a + \xi+ r_-(a)\xi+r_-(a)\xi \ten a+ r_-(a)^2(\xi) \pl .
 \end{align*}
This gives
  \begin{align*}
  &L(s(a)^2)R(s(a)^2-1)(\xi)
   \lel s(a)^2\bigg(\xi\ten a\ten a + r_-(a)\xi+r_-(a)\xi \ten a+ r_-(a)^2\xi\bigg) \\
  &= s(a)\bigg( a\ten \xi \ten a\ten a+ l_-(a)\xi\ten a \ten a
  + (q^n+q^{n+1})\xi\ten a + a\ten r_-(a)\xi + l_-(a)r_-(a)\xi \\
  &\quad  + a \ten r_-(a)\xi\ten a + l_-(a)r_-(a)\xi\ten a
  + q^{n-1} r_-(a)\xi + a \ten r_-(a)^2\xi + l_-(a)r_-(a)^2\xi \bigg) \\
  &= a\ten a \ten \xi \ten a\ten a + \xi\ten a\ten a+
  q a \ten l_-(a)\xi\ten a\ten a + (q^{n+1}+q^{n+2}) a\ten \xi \ten a \\
  &\quad + a \ten l_-(a)\xi\ten a \ten a + a \ten l_-(a)^2\xi \ten a\ten a + (q^n+q^{n+1})l_-(a)\xi\ten a \\
  &\quad + (q^n+q^{n+1})(a\ten \xi \ten a + l_-(a)\xi\ten a+q^n\xi) \\
  &\quad + a\ten a \ten r_-(a)\xi+ r_-(a)\xi+ a\ten l_-(a)r_-(a)\xi \\
  &\quad + a\ten l_-(a)r_-(a)\xi +l_-(a)^2r_-(a)\xi \\
  &\quad   + a\ten a\ten r_-(a)\xi \ten a + r_-(a)\xi\ten a
  + q a \ten l_-(a)r_-(a)\xi\ten a + q^n a\ten r_-(a)\xi \\
  &\quad + a\ten l_-(a)r_-(a)\xi\ten a+ l_-(a)^2r_-(a)\xi \ten a + q^{n-1} l_-(a)r_-(a)\xi \\
  & \quad +  q^n (a \ten r_-(a)\xi + l_-(a)r_-(a)\xi) \\
  & \quad + a \ten a \ten  r_-(a)^2\xi + r_-(a)^2\xi + l_-(a)  r_-(a)^2\xi \\
  &\quad + a\ten l_-(a)r_-(a)^2\xi + l_-(a)^2 r_-(a)^2\xi \pl .
  \end{align*}
The relevant term here is $q^n(q^{n}+q^{n+1})\xi \lel q^{2n}(1+q)\xi$. We get this for $a$ and for $b$ and in $L(x_g^2)R(x_g^2)(\xi)$ as well as $L(x_g^2)R(x_g^2)(\xi)$. As in \cite{Kro} it is easy to prove that
 \[  \max\{ \|\sum_{g\in F} a_g\ten \xi \|_2, \|\sum_{g\in F} \xi\ten a_g \|_2,
 \|\sum_{g\in F} a_g\ten \xi \ten a_g \xi \|_2\}
 \kl C_q\sqrt{|F|} \|\xi\| \]
for \emph{all} $\xi$ in the Fock space.  Using duality and the boundedness of $r_{\pm}(a)$ and $l_{\pm}(a)$, we deduce \[ \|\sum_g L(x_g)R(x_g)+R(x_g)L(x_g)-4(1+q)D(\xi)\|_2
 \kl  C_q \sqrt{|F|} \|\xi\|_2 \]
holds for the diagonal operator $D(\xi)=q^{2n}$ and $\xi \in \cz^{\perp}$. Finally we observe that $1-q^{2n}\gl 1-q^2$ for all $\nen$ and hence
 \[ \|T(\xi)\|\gl 4(1+q)(1-q^2)|F|\|\xi\| -C_q\sqrt{|F|}\|\xi\| \pl .\]
Thus for $|F|$ large enough we find a spectral gap and hence $C$ is a factor. Alternatively, it suffices to assume that $|F|\times \dim(K) \gl k_0(q)$ in order for $C^0$ to be a factor. \qd

\begin{rem}\label{Gzz} For $G=\zz$ (or $G=\zz_m)$ the algebra $\Gamma_q^0(\zz)$ is invariant under $T_t$.
\end{rem}

\begin{proof} Let $g=1$ be the generator of $\zz$.  We observe that $\si(S(e_g))=s(e_g)$ preserves moments. Indeed, $u_g$ commutes with $s(e_g)$ and $\tau(S(e_g)^{\eps_1}\cdots S(e_g)^{\eps_m})$ is $0$ unless we have as many $\eps_k=\emptyset$ and $\eps_k=*$'s. In that case $\tau(u_{g^{\eps_1}} \cdots u_{g^{\eps_m}})=1$. However, in $\Gamma_q(\zz)$ we may write $e_g=\frac{1}{\sqrt{2}}a+ib$ with $a$, $b$ orthogonal and then find that the von Neumann algebra generated by $s(e_g)$ is exactly $\Gamma_q(\ell_2^2(\rz))$ (which is a factor). Note that $\si$ also extends from $\Gamma_q(\{g,g^*\}, K)\to \Gamma_q(\ell_2^2,K)$ for any $K$. Now we compare the Wick words $X_{\si}$ of a monomial $s(e_{g^{\eps_1}})\cdots s(e_{g^{\eps_m}})$ and
$S(e_{g^{\eps_1}})\cdots S(e_{g^{\eps_m}})$. Note that for a pair partition we can only contract $e_{g^{\eps_j}}$ and $e_{g^{\eps_l}}$ if $(\eps_j,\eps_l)\in \{(\emptyset,*),(*,\emptyset)\}$. Then $u_{g^{\eps_j}}u_{g^{\eps_l}}=1$. Since Wick words are $\Gamma_q(\ell_2^2)$, we deduce that same Wick words are in $\si^{-1}(\Gamma_q(\ell_2^2))=\Gamma_q(\{g,g^{-1}\})$. This implies that $\Gamma_q(\{g,g^{-1}\})$ is closed under $T_t$.\qd

\begin{rem} We conjecture that it is enough to assume that
either $|F|\gl 2$ or $|K|\gl 2$ is for $C$ to be a factor, but we don't have a proof at the time of this writing.
\end{rem}

\begin{cor}\label{ug} Let $|F|\times \dim(K)\gl k_0(q)$ and $g\in F$. Then $A\rtimes \Gamma_q(F,K)$ contains unitaries $v_g$ such that $v_gav_g^*=\sigma_{g}(a)$ and $s_q(\delta_g \ten k)=v_gd$ with $d\in C^0$. Moreover the $v_g$'s are orthogonal relative to $A$, i.e. $E_A(v_gv_h^{*})=\delta_{gh^{-1}}$.
\end{cor}

\begin{proof} Let $S_q(g \ten k)=w_g|S_q(g \ten k)|$ be the polar decomposition, for a fixed $0 \neq k \in K$ and $g \in F$. Then  $|S_q(g \ten k)|=(S_q(g \ten k)^*S_q(g \ten k))^{\frac{1}{2}}$ belongs to $C$. Then we note that
 \[ w_ga |S_q(g \ten k)|
 \lel w_g |S_q(g \ten k)|a \lel S_q(g \ten k)a \lel \sigma_g(a)S_q(g \ten k)
 \lel \sigma_g(a)w_g |S_q(g \ten k)| \pl .\]
Thus $w_gaf=\sigma_g(a)af$ holds for the   support of $|S_q(g \ten k)|$. Since $C$ is a factor we may find orthogonal projections $f_j$ and $v_j\in C$ such that $1=\sum_j f_j$, $v_j^*v_j\le f$ and $v_jv_j^*=f_j$. Then we may define $u=\sum_j v_jw_gv_j^*$ which satisfies $uu^*=\sum_j v_jv_j^*=1$ and
 \[ ua \lel \sum_j v_jw_gav_j^*
 \lel \sum_j v_j\sigma_g(a)fv_j^* \lel \sigma_g(a)u \pl .\]
Hence $v_g=u$ is the required unitary. The moreover statement is straightforward and we leave it to the reader.
\qd

\begin{cor} Let $|F|\times \dim(K)\gl k_0(q)$. Then the center of $A \rtimes \Gamma_q(F,K)$ is contained in the center of $A \rtimes [G,G]$.
\end{cor}

\begin{proof} First we want to show that the center of $M$ is contained in $P_0(M)$, the span of the Wick words of order $0$. We consider the $y_g=S(e_g)^*S(e_g)$ and note that for $\xi=W(\eta)u_ha$ we have
 \[ \xi y_g \lel W(\eta)u_h a
 \lel W(\eta)y_g \si_{g^{-1}}\si_g(u_h)a
 \lel W(\eta)y_g u_h a \pl .\]
The same applies to $S(e_g)S(e_g)^*$, and hence our operator $T$ from Lemma \ref{Kro} is a $A \rtimes G$ bimodule map. Hence every element $z$  in the center is understood in the Hilbert module associated with the conditional expectation onto $A \rtimes G$ and hence
$z=T(z)$ implies that $T(z)$ belongs to $P_0(M)$ because $P_0$ is also a $A\rtimes G$ bimodule map. In particular we deduce from Lemma \ref{inf} that $z$ belongs to the center of $A \rtimes [G,G]$. \qd

\begin{rem} We also see that center of $A\rtimes \Gamma_q^0(F,K)$ is contained in  $A \rtimes [G,G]$.
\end{rem}

\begin{cor}\label{factor} Let $|K|$ and $|F|$ as above,  $[G,G]$ be an ICC group or $[G,G]=\{1\}$ (or $q=0$), and assume that the action $\si:G \curvearrowright A$ is ergodic. Then $M=A\rtimes \Gamma_q(F,K)$ is a factor.
\end{cor}

\begin{proof} If $H$ is an ICC group, then the center of $A\rtimes H$ is contained in $A$. Indeed, if $z=\sum_g a_g u_g$ belongs to the center, then we get $z= u_s z u_{s^{-1}} \lel \sum_g \si_s(a_g)u_{sgs^{-1}}$. By comparing coefficients, we deduce that $\|a_h\|_2=\|a_g\|_2$ for every $h=sgs^{-1}$ in the conjugacy class. Since this class is infinite and $a_g\in \ell_2(G;L_2(A))$, we find $a_g=0$ except for $g=1$.

Thanks to the existence of the unitaries $u_g$, we then deduce that $Z(M)\subset A$ is contained in the fixed point subalgebra of $A$. The additional assumption then implies $Z(M)=\cz$.
\qd

\begin{cor}\label{normalizer} Let $M=A \rtimes \Gamma_q(F,K)$ as before. Then $N_M(A)''=M$.
\end{cor}

\begin{proof} It suffices to consider a Wick word $w=W(\xi_{1}\ten \cdots \ten \xi_{m})u_gu_h$, where $\xi_j=\delta_{g_j}\ten k_j$, $g=g_1, \cdots ,g_m$ and $h \in [G,G]$. Since $\{u_g : g \in [G,G] \} \subset N_M(A) \subset M$, we may assume $h=1$. Then we note that
 \[ wa \lel \si_g(a)w \pl .\]
This implies that $|w|\in A'\cap M$, and for the polar decomposition $w=v|w|$ we also see that the support projections $e=vv^*$, $f=v^*v$ are in $A' \cap M$. Let $v_g$ be as in Corollary \ref{ug}. Then we may define the unitary $u=v+(1-e)v_g(1-f)$ which also satisfies $ua=\si_g(a)u$ and hence $u \in N_M(A)$. The assertion follows from $w=u|w| \in N_M(A)A'\subset N_M(A)''$.
\qd

\begin{lemma} Let $A$ be abelian diffuse,  $M=A\rtimes \Gamma_q(G,K)$  with $K$ infinite dimensional separable and assume that the action $G \curvearrowright A$ is free. Then $A'\cap M=A \bar{\ten} N$, where $N \subset M$ is a type $II_1$ factor and in particular $\mathcal Z(A'\cap M)=A$.
\end{lemma}

\begin{proof} Let $v \in \mathcal{U}(A'\cap M)$. Using the Wick space decomposition we may write
 \[ v \lel \sum_{g\in G}(\sum_{m\geq 0} a_g(m)W(\xi_g(m)))u_g.\]
Then, since $av=va$, for all $a \in A$, we have
\[\sum_g (\sum_m aa_g(m)W(\xi_g(m)))u_g \lel \sum_g (\sum_m \sigma_g(a)a_g(m)W(\xi_g(m)))u_g.\]
Hence, for all $g$ and $m$ we must have $aa_g(m)=\sigma_g(a)a_g(m)$, which by freeness is only possible if $g=1$ or $a_g(m)=0$. Thus $v=\sum_m a_1(m)\ten W(\xi_1(m))$, which belongs to the von Neumann subalgebra $A \bar{\ten} N$, where $N$ is spanned by the Wick words of the form $W(\delta_{g_1}\ten k_1 \ldots \delta_{g_m}\ten k_m)$ with $g_1...g_m \in [G,G]$ and $k_1,...,k_m \in K$. Then $N$ is a type $II_1$ subfactor of $M$ (which is in fact isomorphic to $\Gamma_q(\mathcal H)$, for a separable $\mathcal H$).
\qd

\begin{rem} For the remainder of this paper, we will constantly assume that $G$ is an infinite (discrete, countable) group and that $K$ is an infinite dimensional separable Hilbert space, and we will work only with the objects $A \rtimes \Gamma_q^{\pi}(G,K)$, where $\pi$ is either trivial or given by conjugation, unless mentioned otherwise.
\end{rem}

\begin{rem} Denote by $M=A\rtimes \Gamma_q^{\pi}(G,K)$ and by $\tilde{M}=A\rtimes \Gamma_q^{\pi}(G,K\oplus K)$. Using functoriality, we define a 1-parameter group of automorphisms $\al_{\theta}$ of $\tilde{M}$ by
\[\al_{\theta}(aS_q(g,k\oplus k')u_g) = aS_q(g,o_{\theta}(k\oplus k'))u_g, a \in A,k,k' \in K, g \in G,\]
where $o_{\theta}$ is the rotation introduced in 3.1 (6). Note that the automorphism $\beta \in Aut(\tilde{M})$ defined by $\beta(aS_q(g,k\oplus k')u_g)=aS_q(g,k\oplus (-k'))u_g$, satisfies $\beta^2 = id$, $\beta|_M = id_M$, and $\beta \al_{\theta} \beta = \al_{-\theta}$.
\end{rem}

\section{Rigidity}
\begin{defi}{\rm (}see \cite{PoBe}{\rm ,Prop. 4.1)} Let $Q \subset M$ an inclusion of finite von Neumann algebras. We say that the inclusion is rigid, or that $Q$ is relatively rigid in $M$, if for any net $\Phi_{\alpha}:M \to M$ of normal, completely positive, sub-unital, sub-tracial maps such that $||\Phi_{\alpha}(x)-x||_2 \rightarrow 0$ for all $x \in M$, we have that $||\Phi_{\alpha}(x)-x||_2 \rightarrow 0$ uniformly for $x \in (Q)_1$.
\end{defi}

 We consider a similar rigidity property below. First, we pin down notation.   Let $G$ and $H$ be countable discrete groups with $H$ abelian and let $\al : G \to \mathrm{Aut}(H)$ be an action of $G$ on $H$ by automorphisms.   Since $H$ is abelian, its dual $\hat{H}$, the group of characters $\chi: H \to \TT$, is a compact abelian (multiplicative) group with the topology of pointwise convergence.  As $C_{red}^*(H) = C(\hat{H})$, $H$ acts on $C(\hat{H})^*$ by $$\int_{\hat{H}} f(\chi) d(h\mu)(\chi) = \int_{\hat{H}} f(\chi) \langle h, \chi \rangle d\mu(\chi) \text{ for all } h \in H, \mu \in C(\hat{H})^*, f \in C(\hat{H}),$$ The action of $G$ on $H$ induces an action on $\hat{H}$ by $g \cdot \chi := \chi \circ \alpha_{g^{-1}}$, which in turns yields an action on $C(\hat{H})^*$ by $$\int_{\hat{H}} f(\chi)d(\alpha_g\mu)(\chi) = \int_{\hat{H}} f(L_g\chi) d\mu(\chi) \text{ for all } g \in G, \mu \in C(\hat{H})^*, f \in C(\hat{H})$$
 where $(L_g\chi)(h) = \chi(g^{-1}h)$.

 Let us also recall some basic properties from abstract harmonic analysis. For a discrete abelian group $G$ there is a one to one correspondence between positive measures $\mu$ on $\hat{G}$ and positive definition functions  $\phi:G\to \cz$ so that
  \[  \phi(g) \lel \int_{\hat{G}} \chi(g)d\mu(g) \pl .\]
Indeed, if $\mu$ is a positive measure then
 \[ \sum_{g,h} \bar{a}_ga_h \phi(g^{-1}h)
 \lel \int |\sum_h a_h \chi(h)|^2 d\mu(h) \gl 0 \]
shows that the Fourier transform of a measure is positive definite, and the converse is a standard GNS-construction. The same idea also allows us to estimate the distance of two measure by defining
 \[ \|\phi\|_{B(G)} \lel \inf \{\|\xi\| \|\eta\|: \phi(g)=(\xi,\pi(g)\eta)\} \pl .\]
Here the infimum is taken over $^*$-representations $\pi$. It is well-known that $B(G)$ can be identified with the dual $C^*(G)^*$. In particular for an abelian $H$ we find, thanks to amenability,
 \begin{equation}\label{op}
  B(H) \lel C(\hat{H})^* \pl  .
  \end{equation}
With these preliminaries, we will prove:

  \begin{theorem}
     Let $G$ and $H$ be countable discrete groups with $H$ abelian and let $\al : G \to \mathrm{Aut}(H)$ be an action of $G$ on $H$ by automorphisms.   Suppose the following ``strong'' rigidity property, which we call Property (T$^{++}$), holds (see also \cite{IoaRSR}, Prop.4.9 and Thm.5.1): there are finite sets $F_1 \subseteq G$ and $F_2 \subseteq H$ such that for every $\eps > 0$, there is a $\delta > 0$ so that whenever $\varphi \in  C(\hat{H})^*_+$ of norm one satisfies
     \begin{equation}  ||h\varphi - \varphi|| < \delta \text{ for all } h \in F_2 \text{ and } ||\alpha_g \varphi - \varphi|| < \delta \text{ for all } g \in F_1\end{equation}
      then there is a $\hat{\varphi} \in C(\hat{H})^*_+$ of norm 1 such that $$h \hat{\varphi}  = \hat{\varphi} \text{ for all } h \in H \text{ and } ||\varphi - \hat{\varphi}|| < \eps.$$
      Let $M$ be a finite von Neumann algebra containing $A = L(H)$. 
      Assume moreover that for all $g \in G$, there is a unitary $u_g \in \U(M)$ such that $u_gau_g^* = \alpha_g(a)$ for all $a \in A$, where $\alpha:G \curvearrowright A=L(H)$ is the induced action. Then $A \subseteq M$ is rigid.
  \end{theorem}

  \begin{proof}
     Let $\eps > 0$ and $\delta > 0$ be as guaranteed by the rigidity property for groups for $\epsilon' = \frac{\epsilon}{4}$. Put $\delta' := \frac{\delta}{3}$. Let $\Phi: M \to M$ be a subtracial, subunital, ucp map satisfying $||\Phi(x) - x||_2 < \delta$ for all $x \in F_1 \cup F_2$.
      Let $(\mathcal{H}_{\Phi}, \xi_{\Phi})$ be the corresponding pointed Hilbert bimodule (see e.g. \cite{PoBe}, section 1.1).
      Observe that $||u\xi - \xi u|| < \delta$ for all $u \in F_2$ and  $||u_g \xi u_g^* - \xi||_2 < \delta'$ for all $g \in F_1$.

      Define $\hat{\varphi}: H \to \mathbb{C}$ by $h \mapsto (\xi, h\xi h^{-1})$; this is clearly positive-definite and corresponds to $\varphi \in C(\hat{H})^*$ thanks to \eqref{op}. Observe that for all $u \in F_2$ and $h \in H$,
       \begin{align*}
          \hat{u\varphi}(h) = \hat{\varphi}(uh)= (\xi, uh\xi h^{-1}u^{-1}) = (u^{-1}\xi u, h\xi h^{-1}) = (u^{-1}\xi u - \xi, h\xi h^{-1}) + (\xi, h\xi h^{-1}) \pl .
       \end{align*}
    Using the definition of $B(H)$ we deduce that
 \[  ||u\varphi - \varphi||_{C^*(\hat{H})^*} \leq ||u^{-1} \xi u - \xi|| ||\xi|| < \delta' < \delta \pl .\]
     Similarly, for all $g \in F_1$ and all $h \in H$,
      \begin{align*}
         \hat{\varphi \circ \alpha_g}(h) = \varphi(\alpha_g(h)) = (\xi, \alpha_g(h)\xi \alpha_g(h)^{-1}) = (\xi, u_ghu_g^*\xi u_gh^{-1}u_g^*) = (u_g^* \xi u_g, hu_g^* \xi u_gh^{-1}) \\
         = (u_g^*\xi u_g - \xi, hu_g^* \xi u_gh^{-1}) + (\xi, h(u_g^* \xi u_g- \xi)h^{-1}) + (\xi, h \xi h^{-1})
         \end{align*}
      which implies that $|| \varphi  \circ \alpha_g - \varphi||_1 \leq 2||u_g^*\xi u_g - \xi|| < 2\delta' < \delta$. By hypothesis, there is $\tilde{\varphi}$ which is invariant under the action of $H$ by $\alpha$ and satisfies $||\varphi - \tilde{\varphi}||_1 < \eps'$.

      Decompose $\tilde{\varphi} = \psi_n + \psi_s$ into normal and singular parts, with respect to $\varphi$. Then
       $$\epsilon' > ||\varphi - \tilde{\varphi}||_1 = ||\varphi - \psi_n - \psi_s||_1 = ||\varphi - \psi_n||_1 + ||\psi_s||_1$$
     so $||\varphi - \psi_n||_1 < \epsilon'$ and $||\psi_s||_1 < \epsilon'$. Normalize $\psi_n$ to $\psi := \psi_n/||\psi_n(1)||$. Observe that since $|1 - \psi_n(1)| = |\varphi(1) - \psi_n(1)| < \epsilon'$, $$||\varphi - \psi||_1 \leq ||\varphi - \psi_n||_1 + ||\psi_n - \psi||_1 < 2 \epsilon'$$ Moreover, for each $h \in H$,
     $$||\psi \circ \alpha_h - \psi||_1 = ||(\psi - \varphi)\circ \alpha_h||_1 + ||\varphi - \psi||_1 < 4\epsilon'$$

  Because $\psi_n$ is the normal part of $\tilde{\varphi}$, we can view $\eta_\psi := \sqrt{\psi}, \eta_\varphi := \sqrt{\varphi} \in L_2(L(H), \varphi)$. This satisfies
  $||\eta_\psi - \eta_\varphi||_2 \leq ||\psi - \varphi||_1 < \epsilon$ and $||\alpha_h(\eta_\psi) - \eta_\psi||_2 \leq 2||\psi \circ \alpha_h - \psi||_1 < \epsilon$ for all $h \in H$. Define $\pi: L_2(L(H), \varphi) \to \mathcal H$ by $$\pi\left(\sum c_h \lambda(h) \eta_\varphi\right) := \sum c_h h\xi h^{-1}$$
  This is an isometry and satisfies $\pi(\eta_\varphi) = \xi$ and $||\pi(\eta_\psi) - \xi|| = ||\eta_\psi - \eta_\varphi||_2 < \epsilon$. The almost-invariance of $\eta_\psi$ under $H$ by $\alpha$ implies the almost-invariance of $\pi(\eta_\psi)$, as desired.
  \end{proof}

\subsection{Examples}
	Here, we establish that certain classes of groups satisfy our strong rigidity condition (T$^{++}$).   Let $G$ and $H$ be countable discrete groups with $H$ abelian and let $\al : G \to \mathrm{Aut}(H)$ be an action of $G$ on $H$ by automorphisms.  Here, we will view $\hat{H}$ as an additive group of characters $\chi: H \to \RR / \ZZ$. We denote the standard identification $\iota: \RR / \ZZ \to \TT$ (given by $x + \ZZ \mapsto e^{2\pi i x}$).
We will use the $H$-module structure of $\hat{H}$ defined by $(h \cdot \chi)(h') = \chi(hh')$. 

 If there are  finite sets $F_1 \subseteq G$ and $F_2 \subseteq H$ and $\delta > 0$ such that for every $\varphi \in C(\hat{H})^*_+$,
     \begin{equation} \label{rigidT} ||h\varphi - \varphi|| < \delta \text{ for all } h \in F_2 \text{ and } ||\alpha_g \varphi - \varphi|| < \delta \text{ for all } g \in F_1\end{equation}
we say that the action is {\it $(F_1, F_2, \delta)$-invariant.}

 The first class of examples we consider are for a finitely-generated discrete unital commutative ring $R$ with $G = SL_n(R)$, $H = R^n$, and $\alpha$ the natural action of $G$ on $H$ (here $n \geq 2$).  We may identify $\widehat{R^n}$ with $\hat{R}^n$ by $$\langle (\chi_1,\ldots, \chi_n), (r_1, \ldots, r_n) \rangle := \chi_1(r_1) +\chi_2(r_n) + \cdots + \chi_n(r_n)$$
 Under this identification, the dual action of $SL_n(R)$ on $\hat{R}^n$ is given by the adjoint: $$A \cdot \binom{\chi_1}{\chi_2} =  (A^{\text{T}})^{-1} \binom{\chi_1}{\chi_2}$$ where this multiplication is the standard one, defined by the $R$-module structure of $\hat{R}$.

 In this context, $SL_n(R) \curvearrowright R^n$ has {\it property (T$^{++}$) with constant $M$} if there are finite sets $F_1 \subseteq SL_n(R)$ and $F_2 \subseteq R^n$ and $M > 1$ such that for every $\eps > 0$ and  every $(F_1, F_2, \eps/M)$-invariant probability measure $\mu$ on $\hat{R}^n$, there is an $R^n$-invariant probability measure $\nu$ such that $||\mu - \nu|| < \eps$.





  \subsubsection{$SL_2(R) \curvearrowright R^2$} We first consider the case $n = 2$. Our approach here was inspired by \cite{Shalom}. 
   For these group actions, we will consider certain elementary matrices in $SL_2(R)$, of the form $$L(r) :=  \left( \begin{array}{cc} 1 & 0 \\ r & 1 \end{array} \right) ,  \quad U(r) = \left( \begin{array}{cc} 1 & r  \\ 0 & 1 \end{array} \right) \text{ for certain } r \in R.$$
   and the finite set $$F_1(r) := \{L(\pm r), U(\pm r)\} \subseteq SL_2(R).$$
 For this section, we will set
   $$F_2 := \{\pm e_1, \pm e_2\} \subseteq R^2 \text{ where } e_1 = (1, 0), e_2 = (0, 1).$$

  We begin with the group action $SL_2(\ZZ) \curvearrowright \ZZ^2$, where we identify $\hat{\ZZ} = \TT$.
 To prove rigidity, we will need a few lemmas. The first is verbatim from \cite{Shalom}, restated here for completeness. 
     \begin{lemma} \label{Shalom}
       If $\mu$ is a probability measure on the Borel sets of $\RR^2 \setminus \{0\}$, then there is a Borel set $B$ of $\RR^2 \setminus \{0\}$ and $g \in F_1(1)$ such that $|\alpha_g \mu(B) - \mu(B)| \geq 1/4$.
    \end{lemma}

    \begin{lemma} \label{supp}
      If $\mu$ is a $(F_1(1), F_2, \delta)$-invariant probability measure on $\TT^2$ with $\delta \leq 1/20$, then $\mu$ is  supported at 0.
    \end{lemma}

 \begin{proof}
 Assume that $\mu$ is not supported at 0; we argue towards contradiction. Put $Z := (-1/2, 1/2]$ and $X := (-1/4, 1/4)$.  Identify $\TT$ with $Z$ by
    $x \leftrightarrow e^{2\pi i x}.$
  Since
       \begin{align*}
          \delta > ||\pm e_1 \mu - \mu|| = \left| \int_{Z^2} (e^{2\pi i x} -1) d\mu(x,y) \right| & \geq \left| \int_{Z^2} \text{Re}(e^{2\pi i x} - 1) d\mu(x,y) \right| \\ 
          &  \geq \int_{Z \setminus X \times Z} 2\sin^2(\pi x) d\mu(x,y)   \geq \mu(Z \setminus X \times Z)
       \end{align*}
       since $2\sin^2(\pi t) \geq 1$ for all $1/4 \leq |t| \leq 1/2$. Since the same argument holds for $\pm e_2$,
   we must have
     $\mu(X^2) \geq 1 - 2\delta$.

   Define $\mu_1$ on $Z^2$ by $\mu_1(B) := \mu(B \cap X^2)$ for all Borel subsets $B$ of $Z^2$.  By the above we clearly have $||\mu - \mu_1|| \leq 2 \delta$.  For any $g  \in F_1(1)$,
    $$||\alpha_g \mu_1 -\mu_1|| \leq ||\alpha_g(\mu_1 - \mu)|| + ||\alpha_g \mu - \mu|| + ||\mu - \mu_1|| \leq 5\delta$$
    Put $\mu_2 :=\mu_1/\mu(X^2)$, a probability measure on $X^2$.  Then we have
      $$||\alpha_g \mu_2 - \mu_2|| \leq \frac{5\delta}{1 - 2\delta} \leq \frac{5/20}{1- (1/20)^2} < \frac{1}{4}$$
    for each $g \in F_1(1)$. Since we are assuming that $\mu$ (and hence $\mu_2$) has no support at 0 and since $g \cdot X^2 \subseteq Z^2$ for all $g \in F_1(1)$, our measure $\mu_2$ can be viewed as a probability measure on $\RR^2 \setminus \{0\}$ which violates Lemma \ref{Shalom}, a contradiction.
 \end{proof}

\begin{theorem}
  The group action $SL_2(\ZZ) \curvearrowright \ZZ^2$ satisfies  Property (T$^{++}$) with constant 40 and finite sets $F_1(1)$ and $F_2$.
\end{theorem}

\begin{proof}
Let $\eps > 0$ and put $\delta := \eps/40$.  Let $\mu$ be a $(F_1(1), F_2, \delta)$-invariant probability measure on $\TT^2$. Let $\beta := \inf \{\mu(A) : A \text{ Borel}, 0 \in A\}$. Putting $\tilde{\mu} := \frac{\mu - \beta \delta_0}{1 - \beta}$ gives a probability measure on $\TT^2$ satisfying $\mu = \beta \delta_0 + (1 - \beta) \tilde{\mu}$. Since $\delta_0$ is invariant under the action of both $\ZZ^2$ and $SL_2(R)$, for all $g \in F_1(1)$,
$$||\alpha_g(1-\beta)\tilde{\mu} - (1-\beta)\tilde{\mu}|| \leq ||\alpha_g\mu - \mu|| + \beta ||\alpha_g \delta_0 - \delta_0|| < \delta$$
and similarly for $h \in F_2$.
 Thus,
   $$||\alpha_g \tilde{\mu}- \tilde{\mu}|| < (1-\beta)^{-1}\delta \text{ on } F_1(1) \text{ and } ||h \tilde{\mu} - \tilde{\mu}|| <(1-\beta)^{-1}\delta\text{ on } F_2$$ Since $\tilde{\mu}$ is a probability measure on $\TT^2$ with no support at $0$, by Lemma \ref{supp}, $(1-\beta)^{-1}\delta > 1/20$ so $1 - \beta < 20 \delta$ and
   $$||\mu - \delta_0|| = ||(\beta-1) \delta_0 + (1 - \beta)\tilde{\mu}|| \leq 2 (1-\beta) <40 \delta  = \eps$$
as desired.
\end{proof}


Now let $R$ be any finitely-generated, discrete, unital, commutative ring. Here, we prove an intermediate result for polynomial rings that will be used in the sequel. We may identify the dual of $R[t]$ with the group of power series over $\hat{R}$ in $t^{-1}$, which we denote  $R_t := \hat{R}[[t^{-1}]] = \{\bar{\chi} = \sum_{n=0}^\infty \chi_n t^{-n} : \chi_n \in \hat{R}\}$, via $$\left\langle  \sum \chi_n t^{-n}, \sum r_n t^n \right\rangle = \sum_n \langle \chi_n, r_n \rangle $$
 Under this identification, the usual action of $R[t]$ on $\widehat{R[t]}$ becomes
   $$\left( \sum_n r_n t^n \right) \cdot \left( \sum_m \chi_m t^{-m} \right) = \sum_m \sum_{n \geq m} (r_n \cdot \chi_m)t^{n-m}$$

  For generic $\bar{\chi} \in R_t$, $\chi_n$ will denote the coefficient on $t^{-n}$ in the power series corresponding to $\bar{\chi}$.
We will also need $\tilde{R}_t := \hat{R}((t^{-1})) := \{\bar{\chi} = \sum_{n=m}^\infty \chi_n t^{-n} : \chi_n \in \hat{R}, m \in \ZZ\}$, the space of formal ``Laurent series'' over $\hat{R}$. Of course, we may embed $R_t$ in $\tilde{R}_t$. 
 As before, we require a lemmas from \cite{Shalom}, restated here.
\begin{lemma} \label{Shalom2}
  If $\mu$ is a probability measure on $\tilde{R}_t^2 \setminus \{0\}$, then there is a Borel set $B \subseteq \tilde{R}_t^2 \setminus \{0\}$ and $g \in F_1(1) \cup F_1(t) \subseteq SL_2(R[t])$ such that $|\alpha_g\mu(B) - \mu(B)| \geq 1/5$.
\end{lemma}

Our next lemma will give us the ability to induct to obtain property (T$^{++}$) for polynomial rings over $\ZZ$.


\begin{prop} \label{poly}
   Suppose that the action $SL_2(R) \curvearrowright R^2$ has property (T$^{++}$) with constant $M$ and finite sets $F_1 \subseteq SL_2(R)$ and $F_2$.  Then
    \begin{enumerate}
      \item[(a)]  If $\mu$ is a $( F_1 \cup F_1(1) \cup F_1(t), F_2, \delta)$-invariant probability measure on $R_t^2$ for $\delta \leq (20M)^{-1}$, then $\mu$ is supported at 0.
      \item[(b)] The action $SL_2(R[t]) \curvearrowright R[t]$ has property (T$^{++}$) with constant $40M$ and finite sets $F_1 \cup F_1(1) \cup F_1(t)$ and $F_2$.
    \end{enumerate}
\end{prop}

\begin{proof}
(a) Let $\delta \leq (20M)^{-1}$.
 Let $\mu$ be a $(F_1 \cup F_1(1) \cup F_1(t), F_2, \delta/M)$-invariant probability measure on $R_t^2$. Assume that $\mu$ has no support at $0$; we argue toward a contradiction.

    Consider $P : R_t \to \hat{R}$ defined by $P(\bar{\chi}) = \chi_0$. This induces a pushforward measure $\nu = \mu \circ (P^2)^{-1}$ on $\hat{R}^2$. 
Observe
\begin{align*}
  ||\pm e_1 \nu - \nu|| = \left| \int_{\hat{R}^2} (\chi(\pm 1) - 1) d\nu(\chi, \chi')\right| = \left| \int_{R_t^2} (\chi_0(\pm 1) - 1) d\mu(\bar{\chi}, \bar{\chi}') \right| = ||\pm e_1 \mu - \mu|| < \delta/M
\end{align*}
(and similarly for $\pm e_2$) and for $g \in F_1$ and any Borel $B \subseteq \hat{R}^2$,
\begin{align*}
 |\alpha_g \nu(B) - \nu(B)| = |\nu(gB) - \nu(B)| & = |\mu((P^2)^{-1}(gB)) - \mu((P^2)^{-1}(B))| \\ &= |\mu(g(P^2)^{-1}(B)) - \mu((P^2)^{-1}(B))| < \delta/M
\end{align*}
Thus, $\nu$ is $(F_1, F_2, \delta/M)$-invariant, so by hypothesis, $||\nu - \delta_0|| < \delta$.
%

 Put $X := \{\bar{\chi} \in R_t \mid \chi_0 = 0\}$. Then
  $$|\mu(X^2) - 1| = |\nu(\{0\}) - \delta_0(\{0\})| < \delta$$
so $\mu(X^2) > 1 - \delta$.
Define $\mu_1$ on $R_t^2$ by $\mu_1(B) := \mu(B \cap X^2)$ for all Borel subsets $B$ of $R_t^2$. By the above, we clearly have $||\mu - \mu_1|| \leq \delta$. For any $g \in F_1 \cup F_1(1) \cup F_1(t)$,
 $$||\alpha_g\mu_1 - \mu_1|| \leq ||\alpha_g(\mu_1 - \mu)|| + ||\alpha_g\mu - \mu|| + ||\mu - \mu_1|| \leq 3\delta$$
 Put $\mu_2 := \mu_1 / \mu_1(X^2)$, a probability measure on $X^2$. Observe
   $$||\alpha_g \mu_2 - \mu_2|| \leq \frac{3\delta}{1 - \delta} \leq \frac{3/(20M)}{1 - 1/(20M)} = \frac{3}{20M - 1} < \frac{1}{5}$$
 for each $g \in F_1(1) \cup  F_1(t)$. Since we are assuming $\mu$ (and hence $\mu_2$) has no support at 0 and since $gX^2 \subseteq R_t^2$ for all $g \in F_1(1) \cup F_1(t)$, our measure $\mu_2$ can be viewed as a probability measure on $\tilde{R}_t^2 \setminus \{0\}$, which violates Lemma \ref{Shalom2}.

 (b) Let $\eps > 0$ and set $\delta = \eps(40M)^{-1}$.  Let $\mu$ be a $(F_1 \cup F_1(1) \cup F_1(t), F_2, \delta)$-invariant probability measure on $R_t^2$. Let  $\beta := \inf \{ \mu(A): {A \text{ Borel}, 0 \in A}\}$. Putting $\tilde{\mu} := \frac{\mu - \beta \delta_0}{1-\beta}$ gives a probability measure on $R_t^2$ satisfying $\mu = \beta\delta_0 + (1-\beta)\tilde{\mu}$.  Since $\delta_0$ is invariant under the action of both $R^2 $ and $SL_2(R)$, for all $g \in F_1 \cup F_1(1) \cup F_1(t)$,
   $$||\alpha_g(1-\beta) \tilde{\mu}- (1-\beta)\tilde{\mu} || \leq ||\alpha_g \mu - \mu|| + \beta ||\alpha_g \delta_0 - \delta_0|| < \delta$$
   and similarly for $h \in F_2$.
 Thus,
   $$||\alpha_g \tilde{\mu}- \tilde{\mu}|| < (1-\beta)^{-1}\delta \text{ on } F_1 \cup F_1(1) \cup F_1(t) \text{ and } ||h \tilde{\mu} - \tilde{\mu}|| <(1-\beta)^{-1}\delta\text{ on } F_2$$ Since $\tilde{\mu}$ is a probability measure on $R_t^2$ with no support at $0$, by (a), $(1-\beta)^{-1}\delta > (40M)^{-1}$ so $1 - \beta < 40M \delta$ and
   $$||\mu - \delta_0|| = ||(\beta-1) \delta_0 + (1 - \beta)\tilde{\mu}|| \leq 2 (1-\beta) <40M\delta  = \eps$$
as desired.
\end{proof}

\begin{cor} \label{multipoly}
  Let $m \geq 1$. The action $SL_2(\ZZ[t_1,\ldots, t_m]) \curvearrowright \ZZ[t_1, \ldots, t_m]^2$ has Property (T$^{++}$) with constant $40^{m+1}$ and finite sets $F_1(1) \cup F_1(t_1) \cup \cdots F_1(t_m)$ and $F_2$.
  \end{cor}

  \begin{proof}
     Since $\ZZ$ has Property (T$^{++}$) with constant 40 and finite sets $F_1(1), F_2$, by Proposition \ref{poly}, $\ZZ[t_1]$ has Property (T$^{++}$) with constant $40^2$ and finite sets $F_1(1) \cup F_1(t_1)$ and $F_2$. The result follows by induction.
  \end{proof}

We can now show that for any finitely generated commutative unital ring $R$, the group action $SL_2(R) \curvearrowright R^2$ has Property (T$^{++}$).

\begin{theorem} \label{ring}
  Let $R$ be a finitely generated, discrete, unital commutative ring.  Then $SL_2(R) \curvearrowright R^2$ has Property (T$^{++}$) with constant $40^{m+1}$ and finite sets $F_1^2 := F_1(1) \cup F_1(r_1) \cup \cdots \cup F_1(r_m)$ and $F_2$, where $r_0 = 1, r_1, \ldots, r_m$ are the generators of $R$. \end{theorem}

\begin{proof}
  Let $R$ have generators $r_0 = 1, r_1, \ldots, r_m$.  If $m = 0$, then $R = \ZZ$ and the above arguments apply, so assume $m > 0$. 
   Put $R_m := \ZZ[t_1, \ldots, t_m]$.  Consider the surjective ring morphism $\varphi: R_m \to R$ given by $t_k \mapsto r_k$, which induces two ring morphisms, $\varphi^2 : R_m^2 \to R^2$ (defined by $(r,s) \mapsto (\varphi(r), \varphi(s))$) and $\varphi^{(2)}: SL_2(R_m) \to SL_2(R)$ (defined by $\binom{r \ s}{t \ u} \mapsto \binom{\varphi(r) \ \varphi(s)}{\varphi(t) \ \varphi(u)}$) which respect the group action (i.e., $\alpha_{\varphi^{(2)}(g)} \circ \varphi^2 = \varphi^2 \circ \alpha_g$ for all $g \in SL_2(R_m)$). Note that $$\varphi^{(2)}(F_1(1) \cup F_1(t_1) \cup \cdots \cup F_1(t_m)) = F_1(1) \cup F_1(r_1) \cup \cdots \cup F_1(r_m) \text{ and } \varphi^{2}(F_2) = F_2$$
   The map $\varphi$ also induces an injective continuous group morphism $\hat{\varphi}: \hat{R} \to \widehat{R_m}$ via $\hat{\varphi}(\chi) = \chi \circ \varphi$.

  Let $\eps > 0$. Set $\delta := \eps/40^{m+1}$. Suppose $\mu$ is a $(F_1(1) \cup F_1(r_1) \cup \cdots \cup F_1(r_m), F_2, \delta)$-invariant probability measure on $\hat{R}^2$.  Push $\mu$ forward by $\hat{\varphi}^2 : \hat{R}^2 \to \hat{R}_m^2$ to obtain the measure $\nu = \mu \circ (\hat{\varphi}^2)^{-1}$ on $\hat{R}_m^2$.  For all $g \in F_1(1) \cup F_1(t_1) \cup \cdots \cup F_1(t_m)$ and every Borel $B \subseteq \widehat{R_m}^2$,
  \begin{align*}
     |\alpha_g \nu(B) - \nu(B)| = |\alpha_{\varphi^{(2)}(g)}\mu((\hat{\varphi}^2)^{-1}(B)) - \mu((\hat{\varphi}^2)^{-1}(B))| < \delta
  \end{align*}
and
  \begin{align*}
     ||\pm e_1 \nu - \nu||  = \left| \int (\chi_1(\pm 1) - 1) d\nu(\chi_1, \chi_2) \right| & =  \left| \int_{\hat{R}^2} (\chi_1(\varphi(\pm 1)) - 1) d\mu(\chi_1, \chi_2) \right|
     = ||\pm e_1 \mu - \mu||
     \end{align*}
  (similarly for $\pm e_2$).
By Corollary \ref{multipoly}, $||\nu - \delta_0|| < \eps$.  Hence, for any Borel $B \subseteq \hat{R}^2$, by the injectivity of $\hat{\varphi}^2$,
$$|\mu(B) - \delta_0(B)| = |\mu((\hat{\varphi}^2)^{-1}(\hat{\varphi}^2(B))) - \delta_0(B)| = |\nu(\hat{\varphi}^2(B)) - \delta_0(\hat{\varphi}^2(B))| < \eps$$
so  $||\mu - \delta_0|| < \eps$.
  \end{proof}

\subsubsection{Results for $SL_n(R) \curvearrowright R^n$}
We will apply our results above to show that the action $SL_n(R) \curvearrowright R^n$ has Property (T$^{++}$) for  $n \geq 2$.  Throughout this section, we fix a finitely-generated, discrete unital, commutative ring $R$ generated by $r_0 = 1, r_1, \ldots, r_m$ for $m \geq 0$. We also fix
   $$F_1^2 := F_1(1) \cup F_1(r_1) \cup \cdots \cup F_1(r_m) \subseteq R^2$$
and
    $$F_2^n := \{\pm e_k \mid 1 \leq k \leq n\} \subseteq R^n$$
  where $e_k$ is the vector with a 1 in position $k$ and zeros elsewhere.  Note that $F_2^2$ is $F_2$ from the previous section.

  \begin{prop} \label{three}
   The action $SL_3(R) \curvearrowright R^3$ has Property (T$^{++}$) with constant $2\cdot 40^{m+1}$ and finite sets
      $$F_1^3 = \{1 \oplus g \mid g \in F_1^2\} \cup \{g \oplus 1 \mid g \in F_1^2\} \text{ and } F_2^3 = \{\pm e_j \mid 1 \leq j \leq 3\}$$
\end{prop}

\begin{proof}
  Consider $\imath_1, \imath_2: R^2 \to R^3$ defined by $\imath_1(r_1, r_2) = (r_1, r_2, 0)$ and $\imath_2(r_1, r_2) = (0, r_1, r_2)$. These induce $P_1, P_2: \hat{R}^3 \to \hat{R}^2$ defined by $P_1: (\chi_1, \chi_2, \chi_3) \mapsto (\chi_1, \chi_2)$ and $P_2: (\chi_1, \chi_2, \chi_3) \mapsto (\chi_2,\chi_3)$.

  Let $\eps > 0$ and put $\delta = \eps/(2\cdot 40^{m+1})$.  Suppose $\mu$ is a $(F_1^3, F_2^3, \delta)$-invariant probability measure on $\hat{R}^3$. This gives pushforward measures $\nu_j := \mu \circ P_j^{-1}$ on $\hat{R}^2$. For all $g \in F_1^2$ and all Borel $B \subseteq \hat{R}^2$,
  \begin{align*}
  & |\alpha_g\nu_1(B) - \nu_1(B)| 
  = |\alpha_{g \oplus 1}\mu(P_1^{-1}(B)) - \mu(P_1^{-1}(B))| < \delta  \text{ and } \\
  & |\alpha_g\nu_2(B) - \nu_2(B)| = |\alpha_{1 \oplus g}\mu(P_2^{-1}(B)) - \mu(P_2^{-1}(B))| < \delta. \end{align*}For $1 \leq k \leq 2$, $$||\pm e_k \nu_1 - \nu_1|| = ||\pm e_k \mu - \mu|| < \delta \text{ and } ||\pm e_k \nu_2 - \nu_2|| = ||\pm e_{k+1} \mu - \mu|| < \delta.$$  Therefore, by Theorem \ref{ring}, $||\nu_j - \delta_{(0,0)}|| < \eps/2$ for $j = 1,2$. We now show that $||\mu - \delta_{(0,0,0)}|| < \eps$.

Let $B_1 \subseteq \hat{R}^2$ not containing 0. Then
  $$\mu(B_1 \times \hat{R}) = \mu(P_1^{-1}(B_1))  = \nu_1(B_1) = |\nu_1(B_1) - \delta_{(0,0)}(B_1)| < \eps/2$$
and likewise $|\mu(\hat{R} \times B_1)  < \eps/2$. Now let $B \subseteq \hat{R}^3$ be Borel and suppose $(0,0,0) \notin B$. Then $B \subseteq \hat{R}^3 \setminus \{(0,0,0)\} = (\hat{R}^2 \setminus \{(0,0)\} \times \hat{R}) \cup (\hat{R} \times (\hat{R} \times \hat{R} \setminus \{0\}))$ so
  $$|\mu(B) - \delta_{(0,0,0)}(B)| = \mu(B) \leq \mu(\hat{R}^2 \setminus \{(0,0)\} \times \hat{R}) + \mu (\hat{R} \times (\hat{R} \times \hat{R} \setminus \{0\}) < \eps$$
If $(0,0,0) \in B$, then $\mu(\hat{R}^3 \setminus B) < \eps$ so $\mu(B) > 1 - \eps$ and $|\mu(B) - \delta_{(0,0,0)}(B)| = |\mu(B) - 1| < \eps$, as desired.
\end{proof}

%

\begin{prop} \label{induct}
  Let $n \geq 2$ and suppose $SL_n(R) \curvearrowright R^n$ has Property (T$^{++}$) with constant $M$ and finite sets $F_1^n \subseteq SL_n(R)$ and $F_2^n$. Then $R^{n+2}$ has Property (T$^{++}$) with constant $2M$ and finite sets
     $$F_1^{n+2} := \{g \oplus I_2 \mid g \in F_1^n \} \cup \{I_n \oplus g \mid g \in F_1^2\}$$
  and $F_2^{n+2}$.
\end{prop}

\begin{proof}
  View $R^{n+2} = R^n \times R^2$ and $\hat{R}^{n+2} = \hat{R}^n \times \hat{R}^2$. Let $P_1 : \hat{R}^n \times \hat{R}^2 \to \hat{R}^n$ be the projection $(\chi_1, \chi_2) \mapsto \chi_1$ and likewise let $P_2: \hat{R}^n \times \hat{R}^2 \to \hat{R}^2$ be the projection $(\chi_1, \chi_2) \mapsto \chi_2$.

  Let $\eps> 0$ and choose $\delta(2M)^{-1}$. Let $\mu$ be a $(F_1^{n+2}, F_2^{n+2}, \delta)$-invariant probability measure on $\hat{R}^{n+2}$.  The maps $P_j$ induce the pushforward measures $\nu_1 = \mu \circ P_1^{-1}$ on $\hat{R}^n$ and $\nu_2 = \mu \circ P_2^{-1}$ on $\hat{R}^2$. For any $g \in F_1^n$ and any Borel $B \subseteq \hat{R}^n$,
  \begin{align*}
     |\alpha_g\nu_1(B) - \nu_1(B)| 
     & = |\alpha_{g \oplus I_2}\mu(P_1^{-1}(B)) - \mu(P_1^{-1}(B))| < \delta
    \end{align*}
 In the same way, for any $g \in F_1^2$ and any Borel $B \subseteq \hat{R}^2$,
     \begin{align*}
     |\alpha_g\nu_2(B) - \nu_2(B)| &  = |\alpha_{I_n \oplus g}\mu(P_2^{-1}(B)) - \mu(P_2^{-1}(B))| < \delta
    \end{align*}
  For $1 \leq k \leq n$,
  $  ||\pm e_k \nu_1 - \nu_1||  
    = ||\pm e_k \mu - \mu|| < \delta$
  and similarly, for $k=1,2$,
   $||\pm e_k \nu_2 - \nu_2||  
= ||\pm e_{k+n} \mu - \mu|| < \delta$.
Therefore, by hypothesis and Theorem \ref{ring}, $||\nu_1 - \delta_0|| < \eps/2$ and $||\nu_2 - \delta_0|| < \eps/2$.
It remains to show that $||\mu - \delta_{(0,0)}|| < \eps$.

First, let $B_1 \subseteq \hat{R}^n$ and $B_2 \subseteq \hat{R}^2$ be Borel.  If $0 \notin B_1$, then $(0,0) \notin B_1 \times B_2$ so
  $$\mu(B_1 \times B_2) \leq \mu(B_1 \times \hat{R}^2) = \mu(P_1^{-1}(B_1)) = \nu_1(B_1) = |\nu_1(B_1) - \delta_0(B_1)| < \eps$$
  The same holds if $0 \notin B_2$. Now let $B \subseteq \hat{R}^n \times \hat{R}^2$ be Borel and suppose $(0,0) \notin B$. Then $B \subseteq B_0 := (\hat{R}^n \times \hat{R}^2) \setminus  \{(0,0)\} = (\hat{R}^n \setminus \{0\}) \times \hat{R}^2 \cup \hat{R}^n \times (\hat{R}^2 \setminus \{0\})$, so
    $$|\mu(B) - \delta_{(0,0)}(B)| = \mu(B) \leq \mu(B_0) \leq \mu((\hat{R}^n \setminus \{0\}) \times \hat{R}^2 ) + \mu(\hat{R}^n \times (\hat{R}^2 \setminus \{0\})) <  \eps$$
  Conversely, if $B \subseteq \hat{R}^n \times \hat{R}^2$ is a Borel set containing $(0,0)$, then $\mu(B^c) < \eps$ so $\mu(B) > 1 - \eps$ and $|\mu(B) - \delta_{(0,0)}(B)| = |\mu(B) - 1| < \eps$. Therefore, $||\mu - \delta_{(0,0)}|| < \eps$, as desired.
\end{proof}

\begin{theorem} \label{even}
  Let $n \geq 2$.  Then $SL_{n}(R) \curvearrowright R^{n}$ has Property (T$^{++}$). In particular, for $k \geq 1$, $SL_{2k}(R) \curvearrowright R^{2k}$ has Property (T$^{++}$) with constant $2^{k-1} 40^{m+1}$ and finite sets
    $$F_1^{2k} = \{I_{2j} \oplus g \oplus I_{2(k-j-1)} \mid 0 \leq j \leq k-1, g \in F_1^2\} \text{ and } F_2^{2k} = \{ \pm e_j \mid 1 \leq j \leq 2k\}$$
  while $SL_{2k+1}(R) \curvearrowright R^{2k+1}$ has Property (T$^{++}$) with constant $2^k 40^{m+1}$ and finite sets $$F_1^{2k+1} = \{g \oplus I_{2k - 1} , I_{2j+1} \oplus g \oplus I_{2(k- j - 1)} \mid 0 \leq j \leq k-1, g \in F_1^2\} \text{ and } F_2^{2k+1} = \{\pm e_j \mid 1 \leq j \leq 2k +1\}$$
\end{theorem}

\begin{proof}
  Theorem \ref{ring} and Proposition \ref{three} establish the result for $k=1$. Apply Proposition \ref{induct} by induction to finish the proof.
\end{proof}

\begin{cor}For any unital, commutative, finitely generated, discrete ring R, the inclusion
$L(R^n) \subset L(R^n)\rtimes \Gamma_q(SL_n(R),K)$ is rigid.
\end{cor}
\begin{proof}Use theorems 5.2 and 5.12.
\end{proof}
\subsubsection{Related results} We conjecture that property $T^{++}$ holds whenever the inclusion $L(H)\subset L(H\rtimes G)$ is rigid and $H$ is commutative. Indeed, our conditions seems closely related but not obviously identical to the conditions in \cite[Theorem 6.1]{IoaRSR}. We are thankful to Adrian Ioana who brought the results in his paper and the similarities to our approach to our attention. Indeed, in \cite{IoaRSR} Ioana defines an equivalence relation $\mathcal{R}$ on $X$ to be rigid if the inclusion $L^{\infty}(X,\mu)\subset L(\mathcal{R})$ is rigid, and similarly for an action $G$ on $X$, he calls the action rigid if the inclusion $L^{\infty}(X,\mu) \subset L^{\infty}(X,\mu)\rtimes G$ is rigid. In \cite[Theorem 4.4]{IoaRSR} he shows that $G$ acts rigidly on $X$ iff for every sequence of probability measures $\nu_n$ on $X\times X$ such that
 \begin{enumerate}
 \item[i)] $\int f(x)d\nu_n(x)=\int f(x)d\mu(x)=\int f(y)d\nu_n(y)$ for all $f$;
 \item[ii)] $\lim_n f_1(x)f_2(y)d\nu_n(x,y)=\int f_1(x)f_2(x)d\mu(x)$;
 \item[iii)] $\lim_n \|g^*\nu_n-\nu_n\|_{M(X\times X)}=0$ for all $g\in G$;
 \end{enumerate}
one has $\lim_n \nu_n(\{(x,x): x\in X\})\lel 0$.

Let us now consider $M=A\rtimes \Gamma_q(G,K)$ with infinite dimensional $K$ and $T_n$ a sequence of completely positive unital and trace preserving  maps on $M$ such that $\lim_n T_n(x)=x$ holds for all $x\in M$. Since every completely positive map $T_n$ admits a GNS-construction, we find a Hilbert space $H$ and unit vector  $\xi_n\in H$ such that
 \[ \tau(yT_n(x)) \lel \langle \xi_n, x\xi_n y\rangle \pl .\]
In particular
\begin{equation}\label{55}
   \tau(x) \lel \tau(T_n(x)) \lel  \langle \xi_n,x\xi_n\rangle \quad \mbox{and}\quad
\tau(x) \lel \tau(xT_n(1)) \lel  \langle \xi_n,\xi_nx\rangle \pl.
\end{equation}
Let us denote by $\pi_n:M\ten M^{op}\to B(H)$ the corresponding representation given by $\pi_n(x\ten y^{op})(\eta)=x\eta y$. We also assume $A=L^{\infty}(X)$ to be commutative and deduce that the restriction $\pi_n:L^{\infty}(X)\ten_{\min}L^{\infty}(X)\to B(H)$ extends to representation of $L^{\infty}(X\times X)$ and hence there exists a measure $\nu_n$ on $X\times X$ such that
 \[ \int f_1(x)f_2(y) d\nu_n(x,y) \lel \langle \xi_n,f_1\xi_n f_2\rangle \pl .\]
It follows from \eqref{55} and the above that  condition i) is satisfied. The pointwise convergence of $T_n$ to identity for elements $f_1\in A$ implies ii). Let $u_g$ be the unitaries from Corollary 3.14. Since $\lim_n T_n(u_g)=u_g$ we deduce that
 \[ \lim_n \|u_g\xi_n u_g^*-\xi_n\|
 \lel \|u_g\xi_n-\xi_nu_g\| \lel 0 \pl .\]
However, introducing $\xi_n^g=u_g\xi_nu_g^*$ we find that
 \begin{align*}
   (\xi_n^g,\pi_n(f_1\ten f_2)\xi_n^g)
   &= (\xi_n, u_g^*f_1u_g\xi_n u_g^*f_2u_g)
   \lel (\xi_n, \al_g(f_1)\xi_n\al_g(f_2))
   \lel \int \al_g(f_1)\al_g(f_2)d\nu_n \pl .
 \end{align*}
Using the fact that $\nu_n$ is the unique extension from the product $\si$-algebra, we also obtain that for the diagonal action $g^*$ on $X\times X$ we have \[ \|g^*\nu_n-\nu_n\|\kl 2 \|\xi_n^g-\xi_n\| \pl .\]
Hence all the conditions i)to iii) are satisfied. The conclusion easily shows that for the diagonal $\Delta=\{(x,x)|x\in X\}$ we have
  \[ \|u\xi_n-\xi_n u\|^2\lel
  \|\xi_n-u^*\xi_n u\|^2 \kl 4 \nu_n(\Delta^c) \]
for every unitary $u$ in $L^{\infty}(X)$. This is of course uniform convergence. Let us reformulate this result.

\begin{cor}[see \cite{IoaRSR}, Theorem 4.4] Assume that the inclusion  $A=L^{\infty}(X) \subset A\rtimes G$ is rigid and $K$ infinite dimensional. Then the inclusion $A\subset A\rtimes \Gamma_q(G,K)$ is also rigid.
\end{cor}

\begin{rem} {\rm Let $\phi:G_1\to G$ be a surjection and $ L^{\infty}(X) \subset L^{\infty}(X)\rtimes G$ be rigid, then the argument above also shows that $L^{\infty}(X) \subset L^{\infty}(X)\rtimes G_1$ is also rigid, and hence the same is true for our $q$-gaussian algebras. In contrast to the results on the $T^{++}$ property in this section, the proof based on Ioana's work  does not allow us to use a concrete finite set, which could be interesting in some applications (and was our initial goal).
} \end{rem}

\section{Intertwining Results}

\begin{prop} Let $Q \subset M = A \rtimes \Gamma_q^{\pi}(G,K)$ be a von Neumann sub-algebra such that the inclusion $Q \subset M$ is rigid. Let $\tilde{M} = A \rtimes \Gamma_q^{\pi}(G,K \oplus K)$, $(\al_{\theta}) \subset Aut(\tilde{M})$ be the canonical group of automorphisms introduced in Remark 4.21 and assume that either $\mathcal {N}_{\tilde M}(Q)' \cap \tilde M = \mathbb{C}1$ (this is the case if, for example, $Q$ is regular in $M$ since $M' \cap \tilde{M} = \cz 1$) or $Q' \cap \tilde{M} \subset M$. Then there exists a non-zero partial isometry $w \in \tilde M$ such that $wy = \alpha_{\frac{\pi}{2}}(y)w, \forall y \in Q$.
\end{prop}
\begin{proof} In the first case the proof follows verbatim \cite{PoI}, the proof of Theorem 4.1, steps 1 to 3. In the second case one uses the proof of Theorem 4.4, (ii) in the same \cite{PoI}.
\end{proof}

\begin{theorem} Let $Q \subset A\rtimes \Gamma_q^{\pi}(G,K)=M$, $\tilde M = A \rtimes \Gamma_q^{\pi}(G,K \oplus K)$, $\pi: G \rightarrow \mathcal{O}(\ell_{\rz}^2(G))$ being either the trivial representation or the conjugation one. If there exists a non-zero partial isometry $v \in \tilde M$ such that $vy = \alpha_{\frac{\pi}{2}}(y)v, \forall y \in Q$ and one of the following conditions holds
\begin{enumerate}
\item $q=0$;
\item $G$ is a group with the Haagerup property and the inclusion $Q \subset A\rtimes \Gamma_q(G,K)$ is rigid;
\item $\pi$ is trivial, $Q$ is abelian and regular in $M$ and $[G,G]$ is weakly amenable and admits an unbounded 1-cocycle into a mixing non-amenable representation (for the terminology, see e.g. \cite{PoVaI});
\item $\pi$ is trivial, $Q$ is abelian and regular in $M$ and $[G,G]$ is weakly amenable and admits a proper 1-cocycle into a non-amenable representation;
\item $\pi$ is trivial, $Q$ is abelian and regular in $M$ and $[G,G]$ is an weakly amenable, non-amenable, bi-exact group (for the terminology, see \cite{PoVaII}),
\end{enumerate}
then $Q \prec_{M} A$.
\end{theorem}

\begin{proof}We will first prove that $Q \prec_M A \rtimes [G,G]$ in all the situations above. By assumption there exists a non-zero partial isometry in $A \rtimes \Gamma_q(G,K \oplus K)$ such that
 \[  vx \lel \al_{\pi/2}(x)v \]
for $x\in Q$. This implies
 \[ E_{K \oplus 0}(v^*\al_{\pi/2}(x)v)
 \lel E_{K \oplus 0}(v^*v)x \]
for all $x \in Q$. Set $z=E_{K \oplus 0}(v^*v) \neq 0$. Let $\eps\lel \frac{E_{K \oplus 0}(v^*v)}{2}$. Set $B = A \rtimes [G,G]$. According to Lemma \ref{bcontrol}, we can find a finitely generated right $B$-module $H \subset L^2(M)$ such that uniformly for all $x \in (Q)_1$ we have
 \[  \inf_{\xi \in H} \|E_{K \oplus 0}(v^*\al_{\pi/2}(x)v)-\xi\|_2< \eps \pl .\]
Fix a finite orthonormal basis $(\xi_j)_{j=1}^{m}$ of $H$ (see e.g. \cite{PoBe}, 1.4.1). Denoting by $p_H$ the orthogonal projection of $L^2(M)$ onto $H$, it follows that (see \cite{PoBe}, 1.4.1)
\[ p_H(x) = \sum_{j=1}^m \xi_j E_B(\xi_j^{*}x), \|p_H(x)\|_2^{2}=\sum_{j=1}^m \|E_B(\xi_j^{*}x)\|_2^{2}, x \in M. \]
The above inequality implies that uniformly for all $x \in (Q)_1$ we have
\[ ||E_{K \oplus 0}(v^*\al_{\pi/2}(x)v) - p_H(E_{K \oplus 0}(v^*\al_{\pi/2}(x)v))||_2 < \eps \pl. \]
For an arbitrary unitary $u\in Q$ this further implies
 \begin{align*}
 \|z\|_2^{2} &=   \|E_{K \oplus 0}(v^*v)\|_2^{2} \lel  \|E_{K \oplus 0}(v^*v)u\|_2^{2}
  \lel  \|E_{K \oplus 0}(v^*\al_{\pi/2}(u)v)\|_2^{2} \\
  &\le (\|p_H(E_{K \oplus 0}(v^*\al_{\pi/2}(x)v))\|_2 + \eps)^2 \\
  &\le 2(\| \sum_{j=1}^m \xi_j E_B(\xi_j^*E_{K \oplus 0}(v^*\al_{\pi/2}(u)v))\|_2^{2} + \eps^2) \\
  & \lel 2(\|\sum_{j=1}^m \xi_j E_B(\xi_j^*E_{K \oplus 0}(v^*v)u)\|_2^{2} +\eps^2) \\
  &\lel  2\sum_{j=1}^m \|E_B(\xi_j^*zu)\|_2^{2} + \frac{\|z\|_2^{2}}{2}\pl .
  \end{align*}
Therefore we find that for all $u \in \mathcal{U}(Q)$
\[ \frac{\|z\|_2^{2}}{4} \kl \sum_{j=1}^m \|E_B(\xi_j^*zu) \|_2^{2} \pl .\]
By Popa's intertwining criterion (\cite{PoI}, Thm.2.1) this implies $Q \prec_M B=A \rtimes [G,G]$.
Let's now prove that in all of the above situations we have $Q \prec_M A$.
\par {\bf Case 1}. $q=0$. This follows from the argument above and the Lemma 3.8, since in this case we can actually take $H \subset L^2(M)$ to be a finitely generated right $B=A$-module.
\par {\bf Case 2}. Now assume $G$ has the Haagerup property and the inclusion $Q \subset A \rtimes \Gamma_q(G,K)$ is rigid. Let $\phi_n:G \to \cz$ be positive definite functions such that $\lim_{g \rightarrow \infty} \phi_n(g)=0$, $\lim_{n \rightarrow \infty} \phi_n(g)=1$ for all $g \in G$. Denote by $C = \Gamma_q(\ell^2(G) \ten K)$ and note that $A \rtimes \Gamma_q(G,K) \subset (A \bar{\ten} C) \rtimes G$. There exists a non-zero partial isometry $v \in M$ such that $vQv^* \subset vv^*(A \rtimes [G,G])vv^*$. Set $Q_0 = vQv^*$ and $q=vv^*$. Since $Q \subset A \rtimes \Gamma_q(G,K)$ is rigid, it follows that $Q \subset (A \bar{\ten} C) \rtimes G$ is rigid, which further implies that $Q_0 \subset q((A \bar{\ten} C) \rtimes G)q$ is rigid. Define the normal ucp maps $\tilde{\Phi_n}:q((A \bar{\ten} C) \rtimes G)q \to q((A \bar{\ten} C) \rtimes G)q$ given by
\[\tilde{\Phi_n}(q((a \ten c)u_g)q)= \phi_n(g)q(a \ten c)u_g q, a \in A, c \in C, g \in G. \]
It follows that $\tilde{\Phi_n}$ converges uniformly to the identity in the $|| \cdot ||_2$ on the unit ball of $Q_0$. Note that $\tilde{\Phi_n}$ restricts to a sequence of normal ucp maps $\Phi_n$ on $q(A \rtimes [G,G])q$. By restricting to $q(A \rtimes [G,G])q$ and taking into account that $Q_0 \subset q(A \rtimes [G,G])q$, we get that $\Phi_n$ converges uniformly to the identity in the $|| \cdot ||_2$ norm on the unit ball of $Q_0$. As the maps $\Phi_n$ are \emph{compact over A}(see \cite{PoBe}, \cite{OPCartanI}, 2.7), this implies $Q_0 \prec_{A \rtimes [G,G]} A$, by Corollary 2.7 in \cite{OPCartanI}(also see \cite{PoBe}, proof of Thm. 6.2). Hence $Q_0 \prec_M A$ and also $Q \prec_M A$.
\par {\bf Cases 3, 4 and 5}. Assume now that $\pi$ is trivial, $Q$ is abelian and regular in $M$ and $[G,G]$ is either a weakly amenable group admitting an unbounded 1-cocycle into a mixing non-amenable representation or a proper 1-cocycle into a non-amenable representation, or a weakly amenable, non-amenable, bi-exact group. It then follows by \cite{PoI}, Lemma 3.5 that $Q_0$ is regular in $qMq$ and then that $Q_0$ is regular in $q(A \rtimes [G,G])q$. Applying \cite{PoVaI}, Thm. 1.2 in the first two cases and \cite{PoVaII}, Thm. 1.4 in the third one, we see that $Q_0 \prec_{A \rtimes [G,G]} A$, which implies $Q \prec_M A$, as desired.
\qd

\begin{cor}Let $M = A \rtimes \Gamma_q^{\pi}(G,K) = B \rtimes \Gamma_q^{\pi}(G',K')$ with the representation $\pi:G \to \mathcal O (\ell^2_{\rz}(G))$ either trivial or given by conjugation and assume that $A$ and $B$ are abelian, the inclusions $A \subset M$ and $B \subset M$ are rigid, $[G,G]$, $[G',G']$ are ICC groups, and the actions $G \curvearrowright A$, $G' \curvearrowright B$ are free and ergodic. If moreover one of the following conditions holds:
\begin{enumerate}
\item $q=0$;
\item $G,G'$ are groups with the Haagerup property;
\item $\pi$ is trivial, $[G,G]$ and $[G',G']$ are weakly amenable groups which admit unbounded 1-cocycles into mixing non-amenable representations;
\item $\pi$ is trivial, $[G,G]$, $[G',G']$ are weakly amenable groups which admit proper 1-cocycles into non-amenable representations;
\item $\pi$ is trivial, $[G,G]$, $[G',G']$ are weakly amenable, non-amenable bi-exact groups,
\end{enumerate}
then $A \prec_M B$ and $B \prec_M A$ and consequently $\mathcal R(A \subset M)$ and $\mathcal R(B \subset M)$ are stably isomorphic.
\end{cor}
\begin{proof} One applies Proposition 6.1, Theorem 6.2 and Theorem 3.3 in \cite{MeVa}, taking into account that $A$ and $B$ are regular, $\mathcal Z(A'\cap M)=A$, $\mathcal Z(B'\cap M)=B$ by Lemma 4.19 and $M$ is a factor.
\end{proof}

\begin{rem} Proposition 6.2 above and the proof of Theorem 6.3 show that if $A \subset A \rtimes \Gamma_q(G,K)$ is a rigid inclusion, with $[G,G]$ being ICC and $G \curvearrowright A$ ergodic, then whenever $A \rtimes \Gamma_q^{\pi_1}(G,K)=B \rtimes \Gamma_q^{\pi_2}(G',K)=M$, the representation $\pi_2:G' \to \mathcal O(\ell_{\rz}^2(G'))$ being the trivial one, it follows that $A \preceq_M B \rtimes [G',G']$.
\end{rem}

\section{Classification Results}
\begin{prop}Let $M=A\rtimes \Gamma_q^{\pi}(G,K)$. Then $\mathcal R(A \subset M)=\mathcal R(G\curvearrowright A)$.
\end{prop}
\begin{proof} By Corollary 4.14, for every $g \in G$, there exists a unitary $v_g \in M$ such that $v_gav_g^{*}=\sigma_g(a), a \in A$. Take the set $\mathcal F =\lbrace v_g: g \in G \rbrace$. We know that $A'\cap M=A \bar{\ten} N$, where $N$ is generated by the elements $W(\delta_{g_1} \ten k_1 \ten \ldots \ten \delta_{g_m} \ten k_m)$ with $g_1 \ldots g_m \in [G,G]$ and $k_1,...,k_m \in K$. One can easily check that $M=(\mathcal F \cup (A'\cap M))''$ and that the $\|\cdot\|_2$-closed span of $A\mathcal F A$ is isomorphic, as an $A-A$ bimodule, to a direct sum of bimodules of the form $\mathcal H(\sigma_g)$, for $g \in G$. Then according to Prop.3.2 (Lemma 3.4 in \cite{MeVa}) $\mathcal R(A \subset M)$ is generated by the graphs of $\sigma_g$, with $g \in G$. This means $\mathcal R(A \subset M)=\mathcal R(G \curvearrowright A)$.
\end{proof}

\begin{theorem}Assume that $M = A \rtimes \Gamma_q^{\pi}(G,K) = B \rtimes \Gamma_q^{\pi}(G',K')$, the inclusions $A \subset M$, $B \subset M$ are rigid, $[G,G]$ and $[G',G']$ are ICC and the actions $G \curvearrowright A$, $G' \curvearrowright B$ are free and ergodic. Under any of the five sets of conditions in Corollary 6.3, it follows that $\mathcal R(G\curvearrowright A)$ and $\mathcal R(G'\curvearrowright B)$ are stably isomorphic.
\end{theorem}
\begin{proof}We have $\mathcal R(A \subset M)$ is stably isomorphic to $\mathcal R(B \subset M)$ by Corollary 6.3. On the other hand $\mathcal R(A \subset M)=\mathcal R(G \curvearrowright A)$ by proposition 7.1, so the statement follows.
\end{proof}

\begin{cor}
For any non-trivial free product $G=\ast_i G_i$ which satisfies the hypotheses of Corollary 6.3 (in particular if $G$ is weakly amenable or has the Haagerup property), there exist continuously many pairwise non-isomorphic type $II_1$ factors of the form $L^{\infty}(X) \rtimes \Gamma_q^{\pi}(G,K)$. In particular this applies for $G=\mathbb{F}_n, n\geq 2$.
\end{cor}
\begin{proof}Thanks to Thm. 1.3 in \cite{Ga} (see also \cite{PoGa} for the case of free groups) there exist uncountably many (stably) non-OE free ergodic rigid actions $\ast_i G_i \curvearrowright X$. By the same result these actions can be taken such as to coincide on each factor $G_i$ with any prescribed ergodic action of $G_i$. We then use Theorem 7.2 above.
\end{proof}

\begin{cor}Let $G_1,...,G_m$, $H_1,...,H_n$ be ICC groups, each of which either contains a normal non-virtually abelian subgroup with relative property (T) or is a direct product of a non-amenable and an infinite group. Denote by $G=G_1 \ast \ldots \ast G_m$, $H=H_1 \ast \dots \ast H_n$. Assume that $G$ and $H$ are weakly amenable or have the Haagerup property. Let $G \curvearrowright X$, $H \curvearrowright Y$ be two p.m.p. free ergodic rigid actions such that the restriction to each factor is still ergodic. If $L^{\infty}(X)\rtimes \Gamma_q^{\pi}(G,K)$ is isomorphic to $L^{\infty}(Y)\rtimes  \Gamma_q^{\pi}(H,K)$, then $m=n$ and after a permutation of indices we have $\mathcal R(G_i \curvearrowright X)=\mathcal R(H_i \curvearrowright Y)$, for all $i$. In particular, if each of the actions $G_i \curvearrowright X$ is OE-superrigid (e.g. a Bernoulli action of an ICC lattice with property (T) in $\rm Sp(n,1)$), then $G_i \curvearrowright X$ and $H_i \curvearrowright X$ are conjugate, for each $i$.
\end{cor}
\begin{proof}By Thm. 7.2 we have that $\mathcal R(G \curvearrowright X)$ and $\mathcal R(H \curvearrowright Y)$ are stably isomorphic. We then apply Corollary 7.6 in \cite{IPP} and Corollary 6.7 in \cite{CH}.
\end{proof}

\begin{cor}Let $G=\mathbb{F}_{n_1}\times \ldots \times \mathbb{F}_{n_k} \curvearrowright X$, $G'=\mathbb{F}_{m_1}\times \ldots \times \mathbb{F}_{m_l} \curvearrowright X$ be pmp free ergodic rigid actions. If $k\neq l$, then $L^{\infty}(X)\rtimes \Gamma_q^{\pi}(G,K)$ and $L^{\infty}(X)\rtimes \Gamma_q^{\pi}(G',K)$ are non-isomorphic.
\end{cor}
\begin{proof}By Thm.1.16 and 2.12 in \cite{MoSh}, the actions $G \curvearrowright X$ and $G' \curvearrowright X$ are stably non-OE. Note that the statement still holds if one replaces the free groups by torsion-free groups in the class $\mathcal{C}_{\rm reg}$ of Monod and Shalom, as long as they satisfy the assumptions of Corollary 6.3 (for more information about the class $\mathcal{C}_{\rm reg}$, see Sections 1 and 3 of \cite{MoSh}).
\end{proof}

We will end this section with examples where $\mathcal R(G \curvearrowright A)$ almost completely classifies the associated objects.

\begin{prop}\label{back}  Let $\pi$ be the trivial action and $K$ be infinite dimensional. Let $G$ and $G'$ be two groups acting freely on $A$ such that $\mathcal R(G \curvearrowright A)=\mathcal R(G'\curvearrowright B)$. Then $M=A\rtimes\Gamma_0(G,K)$ and $\tilde{M}=A\rtimes\Gamma_0(G',K')$ are isomorphic. \end{prop}

\begin{proof} Let us assume that $A=L^{\infty}(X,\mu)$.
We recall that $[G]=[G']$ means that the equivalence relations $\mathcal{R}\subset X\times X$ for $G$ and $G'$ are the same. Moreover, we have $A\rtimes G=\mathcal{L}(\mathcal{R})$. The map $(t,gt)\to (t,g)$ induces a measure on $\mathcal{R}$ which is the product measure of $\mu$ and the counting measures on the fibers.
Here it is convenient to use the convolution
 \[ F_1\ast F_2(t,s) \lel \sum_{r\sim t} F_1(t,r)F_2(r,s) \pl .\]
Similarly, we may describe the  $(A \ten \Gamma_q(G,K))\rtimes G$ as functions in $L^2(\R,L^2(\Gamma_q(G,K))$, where the convolution above is replaced by pointwise multiplication in $\Gamma_q(G,K)$. Note that it is enough to show that we have the same moments with respect to the generators. The main point is to see that moment
 \[ E_A((S_0(g_1)f_1)\cdots (S_0(g_m)f_m) \]
can be read off the equivalence relation. Here we assume that $f_1,...,f_m$ are projections such that $\al_{g_m}(f_m)=f_{m-1}$, $\al_{g_{m-1}}(f_{m-1})$, ... $\al_{g_2}(f_2)=f_1$, $\al_{g_1}(f_1)=f_m$. Then we find
 \[ E_A((S_0(g_1)f_1)\cdots (S_0(g_m)f_m)
 \lel \tau(s_0(g_1)\cdots s_0(g_m))  f_m
 \lel \sum_{\si \in NC(\{1,...,m\})} \phi_{\si}(g_1,...,g_m) f_m \pl .\]
Using the equivalence model we may understand $S_0(e_g)$ as the function $F(t,s)=\delta_{s,gt}s_0(e_g)$. We note that $\phi_{\si}(g_1,...,g_m)\neq 0$ implies that $g_1\cdots g_m=1$. We recall that by definition  a non-crossing partition contains a neighboring pair. Hence  the only non-zero contributions $\phi_{\si}(g_1,...,g_m)\neq 0$ are obtained by successively eliminating pairs $S_0(g_j)f_jS_0(g_j^{-1})f_{j+1}$, and hence $g_j^{-1}f_{j+1}=f_j$ and $f_{j-1}=f_{j+1}$. Note here that elements in $f_{j+1}$ and $f_j$ are equivalent.

is to have equivalence
equivalence relations of the form $(t,g_1t)$,  $(g_1t,g_2g_1t)$. Thus we are effectively summing over all path using connecting pairs $(f_j,f_k)$ and then obtain a product of projections. Thus the convolution rule
forces us to consider all combinations of pair-wise equivalent projections.
In other words, the expression
  \[ \sum_{\si \in NC(\{1,...,m\})} \phi_{\si}(g_1,...,g_m) \]
counts exactly the number ways $\psi(t,t_1,...,t_m)$
to erase trivial loops (of length $2$) from the string $(t,t_1,....,t_m=t)$ so that eventually we end up with $(t,t)$. Since $\psi(t,t_1,...,t_m)$ only depends on the equivalence relation, we deduce that $M$ and $\tilde{M}$ are isomorphic. \qd

It is unclear whether this result holds for $q \neq 0$. We may, however, construct new classes of examples which only depend in the equivalence relation (if it comes from an action of a discrete group). Let us assume that $R=R_G$ is an equivalence relation on $X$ coming from a group and consider $\tilde{X}=X/R$ the set of representatives. Then we define a $^*$-algebra $\A$  of sections  $f:R\to \bigcup_{\tilde{x}\in X}\Gamma_q(\ell_2(\tilde{x}))$  such that  $f(t,s)\in \Gamma_q(\ell_2([t]))$. Here the product is given by
 \[ f_1\ast f_2(t,s) \lel \sum_{t\sim r} f_1(t,r)f_2(r,s) \pl. \]
The adjoint operation is given by $f^*(t,s)\lel f(s,t)^*$. We may also introduce a trace
 \[ \tau(f) \lel \int \tau_{\Gamma_q}(f(t,t)) d\mu(t) \pl .\]
Indeed, assuming measurability, and $R=R_G$ and that $G$ acts by measure preserving transformations, we have
 \begin{align*}
 \tau(f_1\ast f_2)
 &= \int_X \sum_{r\sim t} \tau(f_1(t,r)f_2(r,t)) d\mu(t) \lel
 \sum_{g}  \int_X \tau(f_1(t,gt)f_2(gt,t))d\mu(t) \\
 &= \sum_{g}  \int_X \tau(f_2(gt,t)f_1(gt,t)d\mu(t)
  \sum_{g}  \int \tau(f_2(t,g^{-1}t)f_1(g^{-1}t,t)d\mu(t) \\
 &=  \sum_{g}  \int \tau(f_2(t,gt)f_1(gt,t)d\mu(t)
 \lel \tau(f_2\ast f_1) \pl .
 \end{align*}
This calculation is correct provided the functions depending on $t$ are measurable. This will be achieved by looking at particular functions in $\Gamma_q(R)$, the $^*$-algebra generated by the canonical embedding of $A=L_{\infty}(X,\mu)$ and  the elements
 \[ S_q(g)(t,s) \lel \delta_{s,g^{-1}t} s_q(e_s) \pl .\]
Note that $s\in [t]$ and hence $S_q(g)\in \A$. We have a natural embedding of $(A\ten \Gamma_q(\ell_2(G)))\rtimes G$ into $\A$. Indeed for $f\in A$ we define $\pi(f)(t,s)=\delta_{t,s} f(t)$. The left-regular representation is given by
 \[ \la_g(t,s) \lel \delta_{s,g^{-1}t} \pl , \]
and finally $\pi(s_q(e_g))(t,s)=\delta_{t,s}s_q(e_{g^{-1}t})$.

\begin{lemma}\label{elem}
\begin{enumerate}
\item[i)] $\la_g\ast \la_h=\la_{gh}$;
\item[ii)] $S_q(g)=\pi(s(e_g))\la_g$;
\item[iii)] $S_q(g)f= \al_g(f)S_q(g)$, $S_q(g)f=\al_{g^{-1}}(f)S_q(g)$;
\item[iv)] $\la_g\pi(s(e_h))\la_g^*=\pi(s_{gh})$, $\la_g\pi(f)\la_g^{-1}=\pi(\al_g(f))$.
\item[v)] Let $\A_G$ be the von Neumann algebras  generated by $S_q(g)$'s for $g\in G$ and $f\in A$. If $\mathcal R(G \curvearrowright A)=\mathcal R(\tilde{G}\curvearrowright A)$, then $A_G=A_{\tilde{G}}$.
\end{enumerate}
\end{lemma}

\begin{proof} Properties i)-iv) are elementary and easy to check. If $[G]=[\tilde{G}]$, then $\phi_g(t)=g^{-1}t$ may be written as
 \[ \phi_g(t) \lel \sum_{h\in \tilde{G}} 1_{E_h}(t) \tilde{\phi}_h(t) \pl ,\]
where  $\tilde{\phi}_h(t)$ uses the action from $\tilde{G}$. This implies that
  \begin{align}\label{ree}
  S_q(g)(t,s) &=  \sum_h 1_{E_h}(t) \delta_{s,h^{-1}t}s_q(e_s)
  \end{align}
is a limit of linear combinations in $\A_{\tilde{G}}$.
\qd

\begin{defi} $\Gamma_q(R_G)$ is the von Neumann algebra obtained from the $GNS$ construction of $A_G$ with respect to $\tau$. Similarly we define $\Gamma_q(R_G,K)$ for an additional real Hilbert space $K$.
\end{defi}

We see immediately that $\Gamma_q(R_G)$, $\Gamma_q(R_G,K)$ are isomorphic to the subalgebra of $A\ten \Gamma_q(\ell_2(G))$ generated by $A$ and the elements
$S_q(g)=s(e_g)u_g$ in the crossed product, and hence these algebras
resemble our previous constructions. The only difference here is that we work with the real Hilbert  space $\ell_2(G)$ instead of the complex version from section 1.

\begin{prop}\label{what} Let $K$ be infinite dimensional, then $\Gamma_q(R_G,K)$ is left invariant by the semigroup $T_t$ given by the number operator. Moreover,
 \[ \Gamma_q(R_G,K) \lel (A \bar{\ten} \Gamma_q(\ell_2(G)\ten K)\rtimes G \pl .\]
\end{prop}

\begin{proof} The proof of the first assertion is the same as in section 1, and we skip it. For the second assertion we fix a sequence  $k_j$ which converges to $0$ weakly and consider the operators
 \begin{align*}
  x_j(t,s)& \lel S_q(g_1\ten k_j)^*S_q(g_2\ten k)S_q(g_3\ten k_j)(t,s) \\
  &= \delta_{r,g_1t} \delta_{v,g_2^{-1}r}\delta_{s,g_3^{-1}v}
   s_q(e_t\ten k_j)s_q(e_v\ten k)s_q(e_{s}\ten k_j) \\
  &=  \delta_{v,g_2^{-1}r}\delta_{s,g_3^{-1}g_2^{-1}gt}
   s_q(e_t\ten k_j)s_q(e_{g_2^{-1}gt}\ten k)s_q(e_{s}\ten k_j) \pl .
   \end{align*}
Passing to the limit we obtain
 \[ \lim_j x_j(t,s) \lel q \delta_{t,s} s_q(e_{g_2^{-1}gt})
 \lel q \pi(s_q(e_{g^{-1}g_2}))(t,s) \pl .\]
Since $q\neq 0$ we see that  $\Gamma_q(R\ten K)$ contains $\pi(\Gamma_q(\ell_2(G)\ten K)$ and $\pi(A)$. Similarly, we consider
 \begin{align*}
  y_j(t,s) &= (S_q(g_1\ten k_j)S_q(g_2\ten k_j)^*)(t,s) \\
  &=   \delta_{r,g_1^{-1}t}s_q(e_r\ten k_j)s_q(e_r\ten k_j)\delta_{s,g_2r} \\
  &\rightarrow_{j\to \infty} \delta_{s,g_2g_1^{-1}t}
  \lel \la_{g_1g_2^{-1}}(t,s) \pl ,
  \end{align*}
weakly in $L_2$.   Hence, we also find the image of $L(G)$, even for $q=0$. But then $\pi(s_q(e_g))=S_q(e_g)\la_{g}^{-1}$ is also in $\Gamma_q(R,K)$ and we find $\Gamma_q(R,K)=(A\ten \Gamma_q(\ell_2(G)\ten K))\rtimes G$ in all cases.
\qd

\begin{theorem}\label{real} Let $A$ be abelian, $|q|<1$ and $K$ infinite dimensional. If $\mathcal R(G\curvearrowright A)=\mathcal R(\tilde{G}\curvearrowright A)$ then $(A \bar{\ten} \Gamma_q(\ell_2(G) \ten K))\rtimes G$ and
$(A \bar{\ten} \Gamma_q(\ell_2(\tilde{G}) \ten K))\rtimes \tilde{G}$ are isomorphic. Conversely, if
 \begin{enumerate}
 \item[i)] $A$ and $\tilde{A}$ are abelian, the inclusions $A \subset (A \bar{\ten} \Gamma_q(\ell_2(G)\ten K))\rtimes G$ and $\tilde{A} \subset (\tilde{A} \bar{\ten} \Gamma_q(\ell_2(\tilde{G}),K))\rtimes G$ are rigid;
 \item[ii)] One of the conditions in Corollary 6.4 holds;
 \item[iii)] $[G,G]$ is ICC and the action of $G$ is free and ergodic,
 \end{enumerate}
then $(A \bar{\ten} \Gamma_q(\ell_2(G)\ten K))\rtimes G \cong (\tilde{A} \bar{\ten} \Gamma_q(\ell_2(\tilde{G})\ten K))\rtimes \tilde{G}$ implies that $\mathcal R(G\curvearrowright A)$ and $\mathcal R(\tilde{G}\curvearrowright \tilde{A})$ are stably isomorphic.

\end{theorem}

\begin{proof} Thanks to Lemma \ref{elem}v), we deduce from $G \subset [\tilde{G}]$ that $S_q(g)$ belongs to the $L_2(A_{\tilde{G}},\tau)$. For $|q|<1$ we also know that $S_q(g)$ is a bounded operator which commutes with the right action and hence $S_q(g) \in \Gamma_q(R_{\tilde{G}})$. The same argument works for $S_q(g \ten k)$. Since $A_G$ is generated by such elements we deduce that $L_2(A_G,\tau)\subset L_2(A_{\tilde{G}},\tau)$ and then
 \[ \Gamma_q(R_G,K) \subset \Gamma_q(R_{\tilde{G}},K)  \pl .\]
Similarly $\tilde{G} \subset [G]$, then implies $L_2(A_G,\tau)=L_2(A_{\tilde{G}},\tau)$ and also
 \[ \Gamma_q(R_{\tilde{G}},K) \subset  \Gamma_q(R_G,K) \pl .\]
Thus we have equality and then Proposition \ref{what} implies the assertion. For the converse, one should note that Theorem 7.2 applies to the objects $(A \bar{\ten} \Gamma_q(\ell^2(G) \ten K)) \rtimes G$, based on verbatim the same arguments as in sections 6 and 7, so the conclusion follows. We leave the details to the reader. \qd

\begin{cor}Let $A$ be abelian, $|q|<1$ and $K$ infinite dimensional. If $\mathcal R(G \curvearrowright A)=\mathcal R(\tilde{G} \curvearrowright A)$ and $\pi:G \to \mathcal{U}(\ell^2(G))$ is  the unitary representation given by conjugation on the complex Hilbert space $\ell_2(G)$, then $A\rtimes \Gamma_q^{\pi}(G,K)$ and $A\rtimes \Gamma_q^{\pi}(\tilde{G},K)$ are isomorphic.
\end{cor}

\begin{proof} Let us first observe that \eqref{ree} implies
\begin{align}\label{eee}
  \pi^{G,\rz}(S_q^{G,\rz}(e_g\ten \xi))
   \lel \sum_{h} \pi(1_{E_h})\pi^{\tilde{G},\rz}(S_q^{\tilde{G},\rz}(e_h\ten \xi))
   \end{align}
for any vector $\xi \in K$. Then we note that the map $v:\ell_2(G;\cz\ten K)\to \ell_2(G;\rz\ten K)+i\ell_2(G;\rz\ten K)\subset \ell_2^{\cz}(G;\rz\ten \ell_2^2\ten K)$ given by
 \[ v(e_g^{\cz}\ten \xi) \lel \frac{e_g\ten e_1\ten \xi + i e_g\ten e_2\ten \xi}{\sqrt{2}} \]
is a real $G$-equivariant map. Indeed, we have
  \[ \la(g) \lel \frac{\la(g)+\la(g)^*}{2} + i \frac{\la(g)-i\la(g)^*}{2i} \]
and    $(wgw^{-1})^{-1}=wg^{-1}w^{-1}$ implies
  \[ \la(w)\la(g)\la(w)^{-1} \lel \al_w(\frac{\la(g)+\la(g)^*}{2}) + \al_w(\frac{\la(g)-i\la(g)^*}{2}) \pl . \]
In particular $A\rtimes \Gamma_q(G,K)$ is canonically embedded in $A\ten \Gamma_q(\ell_2(G;\rz)\ten \ell_2^2\ten K))$ via
 \[ j_{\cz,G}(S_q^{G,\cz}(e_g\ten \xi)) \lel 2^{-1/2} ( S_q^{\rz,G}(e_g\ten e_1\ten \xi)+iS_q^{\rz,G}(e_g\ten e_2\ten \xi)) \pl .\]
Let  $\pi_{G,\tilde{G}}^{\rz}$ be the canonical isomorphism from  Theorem \ref{real}. Then we note that \[ \pi_{G,\tilde{G}}^{\rz}(j_{\cz}(S_q^{G}(e_g\ten \xi))
 \lel \sum_h 1_{E_h}(j_{\cz,\tilde{G}}(S_q^{\tilde{G}}(e_h\ten \xi)) \]
Thus we find
  \[ \pi_{G,\tilde{G}}^{\rz}(A\rtimes \Gamma_q(G,K)) \subset A\rtimes \Gamma_q(\tilde{G},K) \]
and vice versa.  Thus $\tilde{\pi}^{\rz}_{G,\tilde{G}}$ induces indeed an isomorphism between the two algebras.
\qd

\bibliographystyle{amsplain}


\end{document}